\documentclass[12pt]{amsart}

\usepackage{amsthm}

\usepackage{latexsym}
\usepackage{amsmath}
\usepackage{amssymb}
\usepackage{amsfonts}

\theoremstyle{definition}

\newtheorem{dfn}{Definition}[section]

\newtheorem{theorem}[dfn]{Theorem}
\newtheorem{corollary}[dfn]{Corollary}
\newtheorem{lemma}[dfn]{Lemma}
\newtheorem{lem}[dfn]{Lemma}

\newtheorem{proposition}[dfn]{Proposition}
\newtheorem{definition}[dfn]{Definition}
\newtheorem{example}[dfn]{Example}
\newtheorem{rem}[dfn]{Remark}
\newtheorem{question}[dfn]{Question}
\newtheorem{prob}[dfn]{Problem}

\newcommand\Si{\Sigma}
\newcommand\Ga{\Gamma}
\newcommand{\R}{\mathbb{R}}
\newcommand{\Z}{\mathbb{Z}}

\def\<{\langle}
\def\>{\rangle}
\newcommand\al{\alpha}
\newcommand\be{\beta}

\newcommand\La{\Lambda}

\newcommand{\ol}{\overline}

\newcommand{\normal}{\triangleleft}
\newcommand{\acts}{\curvearrowright}

\def\lra{\longrightarrow}
\def\dist{\operatorname{dist}}
\def\embed{\hookrightarrow}

\input diagrams.tex 

\begin{document}

\title{On quasihomomorphisms with noncommutative targets}
\author{Koji Fujiwara and Michael Kapovich}
\date{\today}

\begin{abstract}

We describe structure of quasihomomorphisms from arbitrary groups to discrete groups. We show that all quasihomomorphisms are ``constructible'', i.e., are obtained via certain natural operations from homomorphisms to some groups and quasihomomorphisms to abelian groups. We illustrate this theorem by describing quasihomomorphisms to certain classes of groups. For instance, every unbounded quasihomomorphism to a torsion-free hyperbolic group $H$ is either a homomorphism to a subgroup of $H$ or is a quasihomomorphism to an infinite cyclic subgroup of  $H$.

\end{abstract}

\maketitle

\section{Introduction}

Let $G$ be a group and $H$ be a  group equipped with a proper left-invariant metric $d$ (e.g., a finitely-generated group, equipped with a word metric). A map $f: G\to H$ is called a {\em quasihomomorphism} if there exists a constant $C$ such that
$$
d(f(xy), f(x)f(y))\le C 
$$
for all $x, y\in G$. In the case when $H$ is discrete (and in this paper we limit ourselves only to this class of groups), $f$ is a quasihomomorphism if and only if  the set of {\em defects} of $f$ 
$$
D(f)=\{f(y)^{-1}f(x)^{-1} f(xy): x, y\in G\}
$$
is finite. A quasihomomorphism with values in $\Z$ (or $\R$, equipped with the standard metric) 
is called a {\em quasimorphism}. 

The concept of quasihomomorphisms goes back to S.~Ulam \cite[Chapter 6]{Ulam}, who asked if they are close to 
group homomorphisms. There is a substantial literature on constructing {\em exotic} quasimorphisms, i.e., ones which are not close to homomorphisms, going back to the work of R.~Brooks \cite{Brooks}, see e.g. \cite{C} and references therein; we will refer to quasimorphisms constructed via this procedure as {\em Brooks quasimorphisms}. On the other hand, very little is known about quasihomomorphisms with values in noncommutative groups. 
The first Ulam-stability theorem was proven by Kazhdan \cite{Kazhdan}, namely, that $\epsilon$-quasihomomorphisms from amenable groups into the group of unitary transformations of any Hilbert space are $\epsilon'$-close to homomorphisms (with $\lim_{\epsilon\to 0}\epsilon'=0$). It was proven by Shtern \cite{Shtern} (among other things) that any quasihomomorphism from an amenable group $G$ into $GL(n, \R)$ is a bounded perturbation of a homomorphism. Ozawa \cite{Ozawa} proved that lattices in $SL(n, K)$ ($n\ge 3$, $K$ is a local field) do not admit unbounded quasihomomorphisms to hyperbolic groups. On the negative side, Burger, Ozawa and Thom proved in \cite{BOT} that every group  containing a free nonabelian subgroup, is not {\em Ulam-stable}, in the sense of Kahzdan's paper. Rolli \cite{Rolli} constructed exotic quasihomomorphisms of free groups into groups admitting bi-invariant metrics. After this paper was written, Danny Calegari shared with us an email from Bill Thurston, who noted that``About quasi-morphisms to non-abelian groups: they may be hard to construct in general, but it looks like the Heisenberg group will be one interesting case." In the same email Thurston outlined 
a construction of exotic quasihomomorphisms from hyperbolic 3-manifold groups into the 3-dimensional Heisenberg groups using contact structures on 3-manifolds, although filling in details requires some work; for instance, it is far from clear why quasihomomorphisms defined by Thurston are not close to homomorphisms. 
It follows from our main result that, for this to be the case, at the very least, one has to assume that the 3-manifold $M$ in Thurston's construction satisfies $b_2(M)\ge 2$.  
A construction of quasihomomorphisms (not close to homomorphisms) 
from arbitrary hyperbolic groups to Heisenberg groups, which works  in greater generality, but is purely algebraic and avoids contact structures,  is presented in our Example \ref{ex:nil}. 
 
 Calegari also brought the paper \cite{CaDo} to our attention, where a certain non-commu\-tative version of quasimorphisms into ${\mathbb R}$ is discussed. Furthermore, after this paper was written we received a preprint by Hartnick and Schweitzer \cite{HS}, where they proved existence of exotic quasihomomorphisms of free groups; however, their definition of quasihomomorphisms is different from Ulam's. 
 We will discuss their work in more detail in section \ref{sec:HS}, together with few other generalizations of homomorphisms. In that section we also show that, while Brooks' construction does not generalize to self-quasihomomorphisms of free groups, it does go through when we replace Ulam's notion of a quasihomomorphism with the one of a {\em middle-quasihomomorphism}. 

The goal of this paper is to explain why it is so ``hard to construct" quasihomomorphisms to noncommutative groups which are neither homomorphisms, nor come from quasihomomorphisms with commutative targets, 
provided that $H$ is a discrete group. 

In order to formulate our main theorem we will need a definition:

\begin{definition}\label{def:constructible} 
A quasihomomorphism $f: G\to H$ is {\em constructible} if there exists a finite-index subgroup $G_o<G$, a subgroup $H_o<H$, a finitely generated abelian subgroup $A<H_o$ central in $H_o$,  and a quasihomomorphism $f_o: G_o\to H_o$  within finite distance from $f|G_o$ such that: 

The projection $f': G_o\to Q=H_o/A$ of $f_o$ is a homomorphism. 
\begin{diagram} 
1 & \rTo &   &       &G_o                  & \rTo           &           &     G        &                     &        \\
   &       &    &       &  &\rdTo(2,4)^{f'} &            &                 & &       \\
   &       &    &       & \dDashto^{f_o}                        &               &           &        \dDashto^f         &                     &     \\
      &       &    &       &                         &                 &            &     H          &                      &    \\
1 & \rTo&A & \rTo & H_o                  & \rTo           &  Q      &    \rTo    &  1                  &         \\
\end{diagram}
\end{definition}

\noindent Special subclasses of  quasihomomorphisms include:

1. {\em Almost homomorphisms}, i.e., maps $f: G\to H'< H$, where $H'$ contains a finite normal subgroup 
$K$  such that the projection of $f$ to $H'/K$ is a homomorphism. 

2. Products of quasimorphisms: $f: G\to H'\cong \Z^n < H$; in this case $f=(f_1,\ldots,f_n)$, where each $f_i: G\to \Z$ 
is a quasimorphism. 

\medskip
When we cannot specify the quotient group $Q$ in Definition \ref{def:constructible}, but can only claim that it belongs to a certain class ${\mathcal C}$ of groups, we will say that the quasihomomorphism $f$ in this definition is {\em constructible  from quasihomomorphisms to groups in the class ${\mathcal C}$}.

\medskip
Our main theorem is: 

\begin{theorem}\label{main}
Every quasihomomorphism $f: G\to H$ is {\em constructible}. 
\end{theorem}

We will prove this theorem in section \ref{sec:rigidity} (see Theorem \ref{main1}).

\begin{rem}
Theorem \ref{main} essentially reduces the study of quasihomomorphisms $G\to H$ 
to analyzing quasihomomorphisms $G_o\to A$, homomorphisms $f': G_o\to Q$ and cohomology classes 
$[\omega]\in H^2(Q; A)$ with bounded pull-back classes  $f'^*([\omega])\in H^2(G_o; A)$, 
see section \ref{sec:split}. 
\end{rem}


We also show how one can sharpen the main theorem by restricting to special classes of target groups, e.g.,  some periodic groups (Example \ref{ex:burnside}), hyperbolic groups (Theorem \ref{thm:hyp}), $CAT(0)$ groups (Theorem \ref{thm:cat0}), mapping class groups (Theorem \ref{thm:mcg}) and groups acting on simplicial trees (Lemma \ref{lem:tree}). For instance:

1. All quasihomomorphisms to free Burnside groups $B(n,m)$ (with large odd exponent $m$) 
are bounded. 

2. All unbounded  quasihomomorphisms to hyperbolic groups are either {almost homomorphisms} or have elementary images.  

3. All  quasihomomorphisms $G\to H=Map(\Si)$ to the mapping class group 
are constructible from homomorphisms to other mapping class groups of surfaces (proper subsurfaces in $\Si$), see Theorem \ref{thm:mcg} for the more precise statement. 

In  particular, we will show that higher rank irreducible lattices do not admit unbounded quasihomomorphisms to hyperbolic groups and to mapping class groups. This sharpens the main result of Ozawa in \cite{Ozawa}, since he could prove it only  for lattices in $SL(n,K)$. 

\medskip
Denis Osin \cite{Osin} extended our results on rigidity of quasihomomorphisms to hyperbolic groups and mapping class groups, to the case of relatively hyperbolic target groups and target groups which act acylindrically on Gromov--hyperbolic spaces. Lastly, we note that Nicolaus Heuer in his thesis \cite{Heuer} 
studied quasihomomorphisms to Lie groups. 
 
 \medskip
{\bf Acknowledgements.} We are grateful to the Max-Plank-Institute for Mathematics in Bonn, where this paper was written. We are also grateful to Marc Burger, Danny Calegari, Ursula Hamenst\"adt, Alessandra Iozzi, Gilbert Levitt, Dan Margalit, Igor Mineyev, Denis Osin, Narutaka Ozawa, Andy Putman, Mark Sapir and Andreas Thom for useful conversations, references and corrections. The first author is supported by Grant-in-Aid for Scientific Research (No. 23244005, 15H05739). The second author was also  supported by the NSF grant DMS-12-05312. We are grateful to the referee for useful comments.

\section{Preliminaries}\label{sec:prelim}

In this section we collect some basic  facts about quasihomomorphisms. 

\subsection{Definition and notation}

Throughout the paper (except for \S \ref{sec:HS}), we will be considering quasihomomorphisms to discrete groups, denoted $H$,  equipped with a proper metric $d$ (whose choice will be suppressed in our notation). The reader can think of a finitely-generated group equipped with a word metric as the main example. 
Set $|h|=d(1, h)$. 

\begin{definition}\label{def:ah}
Suppose that a map $f: G\to H$ between groups has the property that $f(G)$ is contained in a subgroup $J<H$, $J$ contains a finite normal subgroup  $K\normal J$, such that the projection $\bar f: G\to \bar J= J/K$ is a homomorphism. We then will refer to $f$ as an {\em almost homomorphism}, it is a homomorphism modulo a finite normal subgroup (in the range of $f$). 
\end{definition}

Clearly, every almost homomorphism is a quasihomomorphism. 

A quasihomomorphism $f: G\to H$ is called {\em bounded} if its image is a bounded (i.e., finite) subset of $H$. 
Note that every map $f: G\to H$ with bounded image is automatically a quasihomomorphism. 

\medskip
A map $f: G\to H$ is a {\em quasiisomorphism} if it is a quasihomomorphism which admits a {\em quasiinverse},  i.e., a quasihomomorphism $f': G\to H$ such that
$$
\dist(f'\circ f, id)< \infty, \quad \dist(f\circ f', id)< \infty. 
$$
Here and in what follows, for maps  $f_1, f_2: X\to Y$ to a metric space $(Y, d_Y)$, 
$$
\dist(f_1, f_2)=\sup_{x\in X} d_Y(f_1(x), f_2(x)). 
$$
A quasiisomorphism is {\em strict} if $f'=f^{-1}$. Two groups $G, H$ are (strictly) quasiisomorphic to each other if there exists a (strict) quasiisomorphism between these groups. 

\medskip 
In what follows we will frequently use the notation ${\mathcal N}_R(S)\subset H$ to denote the $R$-neighborhood of a subset $S$ in a  discrete group $H$ equipped with a proper metric. We will also use the notation $h_1\sim h_2$ for elements $h_1, h_2\in H$ to denote that 
$$
d(h_1, h_2)\le Const
$$ 
where $Const$ is a certain uniform constant (which is not fixed in advance). Instead of the notation $\sim$, we will also write 
write 
$$p\sim_S q$$ 
if $p=qs$ with $p,q \in H, s \in S$ (where the subset $S$ is bounded).
For example, for a quasihomomorphism $f: G\to H$ with $D=D(f)$, by the definition,  
$$
f(ab) \sim_D f(a)f(b)$$ 
for $a,b \in G$.

\medskip
For two quasihomomorphisms $f_i: G_i\to H, i=1, 2$, the notation $f_1\sim f_2$ will mean that the  domain 
of $f_1$ is a finite index subgroup $G_1< G_2$ and that
$$
\dist(f_1, f_2|G_1)<\infty. 
$$

\medskip
For a subset $D$ of a group $H$ and $n\ge 2$ we will  use the notation $D^n$ to denote the subset of $H$ consisting of products of at most $n$ elements of $D$. More generally, for two subsets $A, B\subset H$ we let
$$
A\cdot B=\{ab: a\in A, b\in B\}. 
$$
We will use the notation $D^{-1}$ for the set of inverses of elements of $D$. Then
$$
h\sim_D h'\iff h'\sim_{D^{-1}} h. 
$$

For an element $h\in H$ we let $ad(h)$ denote the inner automorphism of $H$ defined by conjugation via $h$:
$$
ad(h)(x)= hxh^{-1}. 
$$
The map $ad: H\to Inn(H)<Aut(H)$ is a homomorphism; its image $Inn(H)$ is the group of {\em inner automorphisms} of $H$. The quotient group $Out(H)=Aut(H)/Inn(H)$ is the {\em outer automorphism group} of 
$H$. 

Given a group $H$ and its subgroup $A$ we let $N_{H}(A)$ and $Z_H(A)$ denote the normalizer and the centralizer of $A$ in $H$ respectively. For a subgroup $B< H$ we will also use the notation
$$
N_B(A):= N_H(A)\cap B, \quad Z_B(A):= Z_H(A)\cap B. 
$$

\subsection{Elementary properties  of quasihomomorphisms} \label{sec:functoriality}
$~$ 

\medskip 
{\bf Composition.} The composition of quasihomomorphisms is again a quasihomomorphism:
$$
D(f_2\circ f_1)\subset D(f_2)\cdot f_2(D(f_1))\cdot D(f_2). 
$$ 
In particular, if $f_2$ is a homomorphism and $f_2(D(f_1))=\{1\}$, then $f_2\circ f_1$ is a homomorphism. 

\medskip 
{\bf Product construction.} Let $f_i: G\to H_i, i=1,...,n$ be quasihomomorphisms. Then their product
$$
f=(f_1, \ldots , f_n):  G\to H_1\times ... \times H_n
$$
is again a quasihomomorphism. Conversely, given a quasihomomorphism 
$$
f=(f_1,...,f_n): G\to H_1\times ... \times H_n,
$$
in view of the composition property above, each component $f_i$ is again a quasihomomorphism.
 
\medskip 
{\bf Closeness of $f(G)$ and $f(G)^{-1}$.} 
Suppose that  
$$
f: G\to H
$$ 
is a quasihomomorphism. Then for $D=D(f)$ we obtain: 
$$
\epsilon=f(1)=f(1)f(1) s, \quad s\in D
$$
and, hence,
$$
\epsilon=s^{-1}\in D^{-1}. 
$$
Furthermore, for $x\in G$ 
$$
1= f(x x^{-1}) \epsilon^{-1} = f(x) f(x^{-1}) s \epsilon^{-1} , s\in D
$$
which implies that
\begin{equation}\label{eq:inverse}
(f(x))^{-1}=  f(x^{-1}) s', \quad s'\in D^2. 
\end{equation}
In particular, the sets $f(G), (f(G))^{-1}$ are Hausdorff-close to each other.

\subsection{Quasiaction and bounded displacement property} \label{sec:1}

By the definition of a quasihomomorphism, for $D=D(f)$:
$$
f(xyz)\sim_D f(xy) f(z) 
$$
and
$$
f(xyz)\sim_D f(x) f(yz) \sim_{D} f(x) f(y)f(z) . 
$$
In particular,
$$
f(xy) f(z) \sim_{D^{-1}} f(xyz) \sim_{D^2}  f(x) f(y)f(z) 
$$
and, hence,  
\begin{equation}\label{eq:qa}
d(f(xy) h, f(x) f(y) h)\le C_3, \quad \forall h\in f(G), \quad C_3= \max\{|s| : s\in D^2 D^{-1}\}. 
\end{equation}

More precisely, 
\begin{equation}\label{eq:key}
f(xy) h \sim_{D^2 D^{-1}} f(x) f(y) h, \quad h\in f(G), x, y\in G.   
\end{equation}
Therefore, the left multiplication by $f(x)$ defines a {\em quasi-action} of $G$ on $f(G)$. 
The set $f(G)$ is not literally preserved by this quasi-action, but
$$
d(f(x) f(G),  f(G))\le C_1, \quad C_1=\max\{ |s|, s\in D\}, 
$$
for all $x\in G$: For $h=f(y)\in f(G)$, 
$$
f(x) h \sim_{D^{-1}} f(xy)  \in f(G).  
$$
 In view of \eqref{eq:qa}, the defect set $D(f)$ has the property that every element $h\in D(f)$  
quasi-acts on $f(G)$ with bounded displacement.  We define the {\em defect subgroup} $\Delta=\Delta_f$ of $f$ to be the subgroup of $H$ generated by $D(f)$. It is then immediate that every element of $\Delta_f$ (quasi)acts on $f(G)$ with bounded displacement. Equation \eqref{eq:key} shows that there exists a finite subset $D'=D'(f)=D^2 D^{-1}\subset \Delta_f$ such that for every $s\in D=D(f)$, 
\begin{equation}\label{conj}
s h =h s', \quad s'\in D'. 
\end{equation}

\begin{rem}
To verify \eqref{conj}, let $h \in f(G)$ and $s \in D=D(f)$, 
then
$$
h^{-1}sh =f(c)^{-1} f(b)^{-1}f(a)^{-1} f(ab) f(c) \sim _{D^2 D^{-1}}
f(c)^{-1}f(b)^{-1} f(a)^{-1} f(a)f(b)f(c)=1
$$
where $f(c)=h$ and  
$$f(b)^{-1}f(a)^{-1} f(ab)=s.$$
\end{rem}

In particular,
\begin{equation}\label{eq:h-con}
h^{-1} \Delta_f h \subset \Delta_f. 
\end{equation}
Since for every $h\in f(G)$, $h^{-1}\in  f(H) D^2\subset f(H) \Delta_f$ (see equation \eqref{eq:inverse}), 
we conclude that
\begin{equation}\label{eq:h-1-con}
h \Delta_f h^{-1} \subset \Delta_f
\end{equation}
as well. Thus:

\begin{lemma}\label{lem:norm}
The sets $f(G)$ and $f(G)^{-1}$ are contained in $N_H(\Delta_f)$, the normalizer of $\Delta_f$ in $H$. In particular, we obtain a homomorphism 
$$
G\to N_H(\Delta_f)/\Delta_f
$$
whose image is $\langle f(G) \rangle /\Delta_f$. 
\end{lemma}

\medskip 
Let $f: G\to H$ be a quasihomomorphism with the defect subgroup $\Delta_f$. 
As we just proved, the image of $f$ is contained in $N=N_H(\Delta_f)$. 
It follows that  there is no harm in replacing the group $H$ with the  group 
$\<f(G)\>$. We assume from now on that $H=N=\<f(G)\>$; we continue to work with the restriction of the 
original left-invariant metric from the target group of $f$ to $\<f(G)\>$. 

\begin{rem}
We observe that if the group $G$ is finitely-generated, so is the group $\<f(G)\>$: 
It is generated by $f(S)$ and $D(f)$, where $S$
is a finite generating set of $G$. 
\end{rem}

By Lemma \ref{lem:norm}, we also obtain a homomorphism 
$$
\varphi=\varphi_f: G\to Out(\Delta_f)=Aut(\Delta_f)/Inn(\Delta_f) 
$$
given by sending $g\in G$ first to the conjugation automorphism 
$$
\tilde \varphi(g)=ad(f(g))\in Aut(\Delta_f)$$ 
$$
\tilde\varphi(g)(\delta)= f(g) \delta f(g)^{-1}, \quad \delta\in \Delta_f 
$$
and then projecting to the group of outer automorphisms.  (The quasihomomorphism $\tilde\varphi$, of course, in general, is not a homomorphism.) Similarly, by the same lemma, we have an {\em antihomomorphism}
$$
\psi: G\to Out(\Delta_f),
$$
$\psi(g)$ defined by sending $g$ to $\tilde\psi(g)=ad(g^{-1})$ and then  projecting to $Out(\Delta_f)$. In view of \eqref{eq:inverse}, we have
$$
\psi(g)= \varphi(g^{-1}). 
$$

Since $\Delta_f$ is generated by the finite subset 
$D(f)$, the automorphisms $\tilde \varphi(g), \tilde \psi(g)$ are determined by their values on the elements $s\in D(f)$; the images of elements $s\in D(f)$ under $\tilde\varphi(g)$ and $\tilde\psi(g)$ 
belong to a finite subset $D'(f)$ (independent of $g$). Therefore, 
the subset
$$
\tilde\varphi(G)\cup \tilde\psi(G) \subset Aut(\Delta_f)
$$
is finite and, thus, the homomorphism $\varphi$ has finite image. We summarize these simple (but useful) observations in 

\begin{lemma}\label{L0.5}
1. There exists a finite subset $\{y_1,...,y_n\}$ of $H$ such that 
$$
\tilde\varphi(G)\cup \tilde\psi(G)\subset \{ad(y_j): j=1,...,n\}. 
$$
2. There exists a finite-index subgroup 
$G_o<G$ such that $\varphi(G_o)=\{1\}$, i.e., every automorphism $\tilde\varphi(g), \tilde\psi(g)\in Aut(\Delta_f), g\in G_o$ is inner. In particular, we can choose the elements $y_1,...,y_n\in \Delta_f$ such that 
$$
\tilde\varphi(G_o) \cup \tilde\psi(G_o) \subset \{ad(y_j): j=1,...,n\}. 
$$
\end{lemma}


\newpage 
\subsection{Lift and projection}
\subsubsection{Quasi-split exact sequences} \label{sec:split}

Consider an exact sequence
\begin{equation}\label{eq:ses}
1\to A \stackrel{i}{\to} B \stackrel{p}{\to} C\to 1.
\end{equation}
Such a sequence is said to be {\em quasi-split} if there exists a quasihomomorphism $C \stackrel{s}{\to} B$ 
such that  $p\circ s=id$. (More generally, one can allow this composition to  have bounded displacement, but we will not need this.)  In what follows, we will identify $A$ with $i(A)$. 
Given a quasi-splitting $q$ we define the mapping 
$$
q(b)= b^{-1} \cdot \left(s\circ p(b) \right), \quad q: B\to A.   
$$

\begin{lemma}\label{lem:split-q}
If $A$ is central in $B$ then  $q$ is a quasihomomorphism. 
\end{lemma}
\proof Pick $b_1, b_2\in B$ and set $c_i=p(b_i)$,
$$
s(c_i)= a_i b_i, \quad a_i=q(b_i)\in A,  \quad i=1, 2. 
$$
Then 
$$
s(c_1 c_2)= s(c_1) s(c_2) \delta, \quad \delta\in D(s). 
$$
Then,
$$
q(b_1 b_2)= b_2^{-1} b_1^{-1} \cdot s(c_1 c_2)= b_2^{-1} b_1^{-1}  s(c_1) s(c_2) \delta=  
$$
$$
b_2^{-1}  a_1  s(c_2) \delta= a_1  b_2^{-1}   s(c_2) \delta = a_1 a_2 \delta = q(b_1) q(b_2)\delta. \qed 
$$

\medskip 
We continue with the hypothesis of the lemma and define the maps 
$$
F: B\to C\times A, \quad F(b)= (p(b), q(b))  
$$
and
$$
F': C\times A\to B, \quad F'(c,a)= s(c) a^{-1}. 
$$
Since $p$ and $q$ are (quasi)homomorphisms, so is $F$. The proof that $F'$ is a quasihomomorphism is completely analogous to the proof of Lemma \ref{lem:split-q} and is left to the reader. 

\begin{lemma}
 $F' =F^{-1}$; in particular, the group $B$ is strictly 
quasiisomorphic to $C\times A$. 
\end{lemma}
\proof $F'\circ F(b)= F'(p(b), q(b))= sp(b) \cdot (q(b))^{-1}=  sp(b) \cdot sp(b)^{-1} \cdot b=b$. The reader will verify that $F\circ F'=id$. \qed 

\medskip 
Given a quasi-split extension \eqref{eq:ses}, 
each quasihomomorphism $f: G\to C$ {\em lifts} to a quasihomomorphism 
$\tilde f: G\to B, \quad \tilde f= s\circ f$.
\begin{diagram}
G & \rDashto^{\tilde{f}} & B    \\
    &  \rdDashto^f          &\dTo^p  \\
    &                               & C\\
\end{diagram}

Similarly, given a quasi-split exact sequence \eqref{eq:ses}, each quasihomomorphism $f: B\to H$ 
{\em projects} to a quasihomomorphism $\bar f=f\circ s: C\to H$.
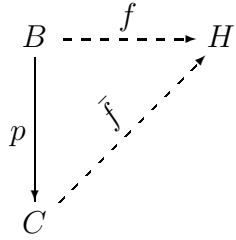
\begin{diagram}
B & \rDashto^{f} & H    \\
\dTo^p    &  \ruDashto^{\bar{f}}          & \\
C    &                               & \\
\end{diagram}

If $f: G\to C$ is unbounded, the quasihomomorphism $\tilde f$ is, of course, unbounded as well. This is not necessarily the case for projections of quasihomomorphisms $C \stackrel{f}{\to} H$ as one can take, for instance, 
$B=A\times C$ and $f=f_1\times f_2: G\to B$, with bounded $f_2$ and unbounded $f_1$. However, if $A$ is finite and $f$ is unbounded, then $\bar f$ is unbounded as well. We will use this observation several times in the case when $H=\Z$, in order to construct unbounded quasimorphisms on the quotient group $C$. 

\medskip
\begin{example}\label{ex:qsplit}
Examples of quasi-split sequences  are given by:

 a. Extensions with finite kernel $A$: In this case {\em any} section $s: C\to B$ will define a quasi-splitting.

b. Central extensions whose obstruction class is a bounded 2nd cohomology class, cf. \cite{Gersten} or \cite{NR}. 

\end{example}
To justify (b), suppose that $\omega\in Z^2(C,A)$ is a bounded {\em normalized} cocycle, i.e.,  $s(1,c)=s(c,1)=0\in A$ for all $c\in C$. Here and in what follows we use the restriction of the metric from $B$ to $i(A)\cong A$. We also refer the reader to \cite{C} for the discussion of bounded cohomology. 

Following \cite[p. 92]{Brown}, we define the extension $E_\omega$ of $C$ by $A$, using the  group law on the product $A\times C$ given by the formula:
$$
(a_1, c_1)(a_2,c_2)= (a_1+a_2+ \omega(c_1,c_2), c_1c_2).
$$
The group $E_\omega$ is then a central extension of $C$ by $A$, which is isomorphic to the one in \eqref{eq:ses}. 
The quasi-splitting of the sequence
$$
0\to A\to E_\omega\to C\to 1
$$
is given by $s(c)=(0,c)$. Then $\omega$ is bounded if and only if $s$ is a quasihomomorphism. We obtain

\begin{lemma}
The sequence \eqref{eq:ses} quasi-splits if and only if the extension class is bounded. 
\end{lemma}

In \S \ref{sec:mcg} we will prove Proposition \ref{prop:map-split} about quasi-splitting of a central extension associated with a certain subgroup of the mapping class group of a surface, illustrating this result.

\subsubsection{Second bounded cohomology of $G$}

Note that there are situations when the sequence \eqref{eq:ses} does not quasi-split, but {\em homomorphisms}  
$f: G\to C$ still lift to quasihomomorphisms $\tilde f: G\to B$. Namely, assume that the subgroup 
$i(A)$ is central in $B$ and the class $f^*([\omega])\in H^2(G; A)$ is bounded. Then the homomorphism $f$ lifts to a quasihomomorphisms $\tilde f: G\to B$. To see this, consider the central extension of $G$ by $A$ defined by the class $f^*([\omega])$:
$$
0\to A \to \tilde E \to G\to 1. 
$$
Let $\tilde s: G\to \tilde E$ be the quasi-splitting. Composing $\tilde s$ with the natural homomorphism 
$\hat{f}: \tilde E\to B$ (which projects to $f: G\to C$), we obtain the required lift $\tilde f$.  The converse to this is also easy to see: 
If $f$ lifts to a quasihomomorphism $\tilde f$, then the class $f^*([\omega])\in H^2(G; A)$ is bounded. 

\begin{example}
Consider the case where $A$ is a finitely-generated abelian group and the group 
$G$ is hyperbolic. Then all cohomology classes in $H^2(G; A)$ are bounded (see \cite{NR}), which implies 
that quasihomomorphisms $f: G\to C$ always lift to quasihomomorphisms $G\to B$. 
\end{example}

\begin{example}\label{ex:nil}
Consider the integer Heisenberg group $B=H_{2n}$, where $A\cong \Z$, $C\cong \Z^{2n}$ and the obstruction class $[\omega]$ is unbounded (the cocycle $\omega$  is the restriction of a symplectic form from $\R^{2n}$ to $\Z^{2n}$). Then  every homomorphism $f: G\to \Z^{2n}$ from a hyperbolic group $G$,  lifts to a quasi-homomorphism $\tilde f: G\to H_{2n}$. We now explain how to use this in order to construct 
examples of quasihomomorphisms to nilpotent groups which are not close to homomorphisms.  

It follows from the definition of $H_{2n}$ that two elements $b, b'\in B$ commute if and only if $\omega(p(b), p(b'))=0$. Take $G$ which admits an epimorphism $f: G\to C'\cong \Z^{2} < \Z^{2n}$ such that $\omega$  is nondegenerate on $C'$ and $f^*(\omega)$ defines a trivial cohomology class of $G$. For instance, 
we can take $G$ to be the fundamental group of a closed oriented surface of genus $\ge 2$ and 
$f: G\to C$ induced by a map of nonzero degree $S\to T^2$. Or, in line with Thurston's suggestion mentioned in the introduction, we can take $G$ to be the fundamental group of a closed hyperbolic 3-manifold $M$ which admits a retraction $r: M\to S$ to a closed oriented hyperbolic surface $S\subset M$. (It follows from the work of 
Agol, Haglund and Wise that for every quasifuchsian surface subgroup of $\pi_1(S)< \pi_1(M)$ there exists a finite index subgroup of $\Ga'< \pi_1(M)$ which retracts to $\pi_1(S)\cap \Ga'$. Hence, examples which we need abound.) Then take the composition of $r$ with a homomorphism induced by a nonzero degree map $S\to T$. 
\end{example}

\begin{lemma}
Suppose that $G$ is a hyperbolic group, $f: G\to C$ is a homomorphism such that $[f^*(\omega)]\ne 0$ in 
$H^2(G, \Z)$. Then:

1.  For each  quasihomomorphism $\tilde f: G\to B$ as above, there is no finite index subgroup $G_o< G$ such that $\tilde{f}|G_o$ is within finite distance from a homomorphism. 

2. The image of $\tilde f$ is not Hausdorff-close to an abelian subgroup of $B$. 
\end{lemma}
\proof 1. Suppose, for the sake of a contradiction, that there exists such $G_o< G$ and a homomorphism 
$f': G_o\to B$ within finite distance from $f|G_o$. Then the distance between the homomorphisms 
$f_o:=p\circ f'$ and $f|G_o$ is again bounded, which implies (since $C$ is free abelian of finite rank) 
that the two homomorphisms are actually equal. Since $G_o$ has finite index in $G$, the transfer argument shows that $[f_o^*(\omega)]=[f^*(\omega)]\in H^2(G_o, \Z)$ is still nonzero. However, for arbitrary central extension 
$$
1\to A\to \tilde{\Gamma} \to \Gamma \to 1 
$$
and arbitrary group $\Lambda$ we have that a homomorphism $h: \Lambda\to \Gamma$ lifts to a homomorphism $\tilde h: \Lambda\to \tilde\Gamma$ if and only if the pull-back $h^*(\omega)$ of the extension cocycle, vanishes in $H^2(\Lambda, A)$. Thus, in our situation, we obtain a contradiction with the assumption about nontriviality of the $f^*(\omega)$. 

2. Suppose that $\tilde{f}(G)$ is Hausdorff-close to an abelian subgroup $B'< B$. Then the subgroup $f(G)< C$ is Hausdorff-close to the abelian subgroup $C'=p(B')$. Since subgroups of the abelian group $C$ are Hausdorff-close iff they are commensurable, we can assume, after replacing $G$ with a  finite index subgroup $G_o<G$, 
that $f(G_o)$ is contained in $C'$ and, hence, $\tilde{f}(G_o)$ is contained in $B'$. As in Part 1, the restriction of the extension class $\omega$ to the finite index subgroup $C_o:=f(G_o)< f(G)$ is still nontrivial.  This, however, implies that each abelian subgroup of $p^{-1}(C_o)$, such as $B'\cap p^{-1}(C_o)$, projects to a cyclic subgroup of $C$, in particular, the restriction of $\omega$ to $p(B')=C_o$ is trivial in this case.  A contradiction. \qed

\begin{rem}
As a warning to the reader, we note that, in general, even if $B$ is finitely-presented, its center  
 may fail to be finitely generated, see e.g. \cite{Abels}. 
\end{rem}

\begin{question}
Is it true that for arbitrary (countable) abelian group $A$ and a hyperbolic group $G$, every class in $H^2(G; A)$ is bounded (i.e.,  is represented by a cocycle taking only finitely many values)? 
\end{question}

Suppose that $\tilde f_1, \tilde f_2: G\to B$ are distinct quasihomomorphisms lifting $f: G\to C$. Then for every 
$g\in G$
$$
\tilde f_2(g) = \phi(g) \tilde f_1(g), 
$$
where $\phi(g)\in A$ (which we identify with $i(A)$). It is immediate that $\phi: G\to A$ is a quasihomomorphism. We summarize these observation as 

\begin{lemma}
Assume that $i(A)$ is central in $B$.  Then: 

1. A homomorphism $f: G\to C$ lifts to a quasihomomorphism $\tilde f: G\to B$ if and only if the pull-back class $f^*([\omega])\in H^2(G; A)$ is bounded. 

2. Different quasihomomorphic lifts differ by quasihomomorphisms $G\to A$. 
\end{lemma}

\subsection{Summary of constructions of  quasihomomorphisms}

So far, we saw several basic constructions of quasihomomorphisms: 

\medskip 
i) {\bf Lift.} If $\bar f: G\to \bar H$ is a quasihomomorphism and $1\to K\to H\to \bar H\to 1$ is a short exact sequence with a (virtually) abelian group $K$, then lift $\bar f$ (if possible) to a quasihomomorphism $f: G\to H$. Note that if the exact sequence quasi-splits with a quasi-splitting $s: \bar H\to H$, then we can always lift $\bar f$ to a quasihomomorphism $f=s\circ \bar f$. For instance, all almost homomorphisms $G\to H$ appear in this fashion. 

ii) {\bf Product.} If $f_i: G\to H_i$ are quasihomomorphisms, $i=1,...,n$, then take
$$
f=(f_1,\ldots,f_n): G\to H= \prod_{i=1}^n H_i. 
$$

iii) {\bf Composition.} The special case of the composition construction is when $f: G\to H$ is a quasihomomorphism and $\iota: H\to \tilde H$ is a monomorphism; then we extend $f$ to the  quasihomomorphism $\tilde f=\iota\circ f$. 

iv) {\bf Extension from a finite index subgroup.} Extend $f_o: G_o\to H$ (if possible) to a quasihomomorphism $f: G\to H$, where $|G:G_o|<\infty$.  

v) {\bf Bounded perturbation.} Replace $f$ (if possible) with a quasihomomorphism $f'$ within finite distance from $f$. Note, however, that (unlike quasimorphisms to abelian groups) a bounded perturbation of a quasihomomorphism need not be a quasihomomorphism. For instance, we will show in Theorem \ref{thm:def-rig} that if $f_1, f_2: G\to H$ are quasihomomorphisms to a torsion-free hyperbolic group,  and $\dist(f_1, f_2)<\infty$, then either $f_1=f_2$, or $f_1, f_2$ are both bounded, or both are quasimorphisms to the same cyclic subgroup. 
Nevertheless, we will see and use repeatedly in the paper that sometimes quasihomomorphisms can be perturbed to quasihomomorphisms. 

\medskip 
By using repeatedly these constructions one can obtain new quasihomomorphisms  from a given set of quasihomomorphisms. In Theorem \ref{main} we show that {\em all quasihomomorphisms are constructible}; in particular, there is no need to repeat the above constructions.  Another construction which, as it turns out, to be not needed (in full generality) is the composition of quasihomomorphisms. One needs only its special cases as in (i) and (iii).

\section{Rigidity of quasihomomorphisms}\label{sec:rigidity}

\subsection{Quasihomomorphisms and centralizers} 

Consider a quasihomomorphism $f: G\to H$. 
By Part 1 of Lemma \ref{L0.5}, there exists a finite subset 
$\{y_1,\ldots, y_n\}\subset G'=f(G)\subset H$, such that 
for every $x\in G$ there exists $y_j$ for which 
$$
\tilde\psi(x)= ad(y_j)\in Aut(\Delta_f),  
$$
i.e., for every $\delta\in \Delta_f$,
$$
f(x)^{-1} \delta f(x)= y_j \delta y_j^{-1}, 
$$
and, hence, 
$$
[f(x)y_j, \delta]=1. 
$$
In other words, $f(x)y_j$ belongs to $Z_H(\Delta_f)$, the centralizer of $\Delta_f$ in $H$. 
Moreover, by Part 2 of the same lemma, if $\varphi(x)=1$ then we can choose 
$y_j\in \Delta_f$.   Recall that the image of the homomorphism $\varphi$ is finite and the kernel 
$G_o=\ker(\varphi)$ has finite index in $G$. 

We, thus, obtain the following strengthening of Lemma \ref{L0.5}:

\begin{corollary}\label{cor:cent}
There exists a constant $C$ such that
$$
f(G)\subset {\mathcal N}_C(Z_H(\Delta_f)). 
$$
Moreover, setting $G_o=\ker(\varphi)$, we get  
$$
f(G_o)\subset \bigcup_{i=1}^n  Z_H(\Delta_f) \cdot y_i, \quad y_i\in \Delta_f.  
$$
\end{corollary}

In particular, 

\begin{corollary}
Suppose that $H$ has the property that the centralizer of every nontrivial element is abelian. Then for every quasihomomorphism $f: G\to H$ either $f$ is a homomorphism or its image lies in a $C$-neighborhood of some abelian subgroup (with $C$ depending on $f$, of course). 
\end{corollary}

\begin{example}\label{ex:burnside}
Let $H$  be either an (infinite) free Burnside group $B(n,m)$ on $n$ generators and odd exponent $m\ge 665$, or a Tarski monster (see \cite{Olshansky}), where all proper subgroups are finite cyclic. 
Note that by a theorem of Adyan and Novikov (see e.g. 
\cite{Olshansky}), the centralizer of every nontrivial element of $B(n,m)$ is cyclic of order $m$.  
In the case of Tarski monsters constructed by Olshansky, centralizers of nontrivial elements are again cyclic,  
Theorem 26.5 of \cite{Olshansky} (we owe the reference to Denis Osin). Therefore, for every group $G$, every unbounded quasihomomorphism $f: G\to H$ is  a homomorphism. (Since if $D(f)\ne \{1\}$ then $f(G)$ is close to the centralizer of $D(f)$.)  
\end{example}

Note, however, that for $m$ even, some centralizers in $B(n,m)$ are infinite,  
see \cite{IO} for the details. This leads to

\begin{question}
Are there quasihomomorphisms $f: G\to H$ to torsion groups $H$, 
which are not within finite distance from almost homomorphisms? 
\end{question}

We note that if $H$ is a nilpotent torsion group, then indeed, the answer to this question is negative (since the defect subgroup is finite in this case). Furthermore, by repeating the construction in Example \ref{ex:nil} with 
$A=\Z_2$ and $G$ a countably infinite direct sum of $\Z_2$'s, it is easy to construct examples of quasi-homomorphisms to torsion nilpotent groups which are not close to homomorphisms. 

\medskip
We next explain how one can alter $f$ such that its image is actually contained in $Z_H(\Delta_f)$. 
As above, let $G_o=\ker(\varphi)$. We let $r: f(G_o)\to Z_H(\Delta_f)$ be a nearest-point projection and set
$$
f_o:= r\circ f: G_o\to   Z_H(\Delta_f) < Z_H(\Delta_{f|G_o})
$$
Clearly, $d(f, f_o)=R<\infty$. 

\begin{lemma}\label{L0}
The map $f_o$ is a quasihomomorphism and $D(f_o)\subset \Delta_f$.  
\end{lemma}
\proof We have
$$
f(x_1 x_2)= f(x_1) f(x_2) s, \quad s\in D(f)
$$
$$
f(x_i)= f_o(x_i) \delta_i, \delta_i\in \Delta_f, f(x_1 x_2)= f_o(x_1 x_2) \delta_3, 
\quad |\delta_i|\le R, i=1, 2, 3. 
$$
Since $f_o(x_i)$ commutes with $\Delta_f$, 
$$
f_o(x_1) f_o(x_2) \delta_1 \delta_2 = f_o(x_1) \delta_1 f_o(x_2)\delta_2= 
$$
$$
f(x_1) f(x_2)= f(x_1 x_2) s= f_o(x_1x_2) \delta_3 s.
$$
Therefore,
$$
f_o(x_1) f_o(x_2) \sim_{D_o} f_o(x_1x_2) 
$$
where $D_o=D(f_o)\subset \Delta_f$ is finite (since $|\delta_i| \le R$
and $s \in D(f)$). \qed 

\medskip
We can now prove

\begin{theorem}\label{main1}
Every quasihomomorphism $f: G\to H$ is constructible:  For the subgroup $G_o< G$ and the 
quasihomomorphism 
$$
f_o: G_o\to H_o<Z_H(\Delta_{f_o})<H$$ 
as above, we have:  

a) The projection of $f_o$ to $\bar f_o: G_o\to  Q=H_o/\Delta_{f_o}$ is a homomorphism. 

b) $H_o=\langle f_o(G_o)\rangle$ and the finitely generated abelian subgroup $\Delta_{f_o}$ is central in $H_o$. 
\end{theorem}
\proof Let $H_o<H$ be the subgroup generated by 
$f_o(G_o)$. By the construction,
$$
f_o(G_o)\subset Z_H(\Delta_f) < Z_H(\Delta_{f_o})
$$
since $\Delta_f > \Delta_{f_o}$. Since $H_o=\<f_o(G_o)\>$, the subgroup $\Delta_{f_o}<H_o$ is central in $H_o$. Since $\Delta_{f_o}$ contains the defect set of $f_o$, the map $\bar f_o$ is a homomorphism. \qed 

\medskip
We note that Theorem \ref{main} from the introduction follows immediately. 

\subsection{Quasihomomorphisms close to abelian subgroups} 

In this and the following section we establish two technical results, which are 
variations of Theorem \ref{main} and will be used in the proof of Theorem \ref{thm:mcg}. 

Let $B$ be a group which is an extension
$$
1\to A \to B \stackrel{p}{\to} C\to 1, 
$$
where $A$ is a finitely generated abelian group. Suppose, further, that $A$ is {\em virtually central} in $B$ in the sense that there exists a finite index subgroup $C'\normal C$ which acts trivially on $A$. We will then refer to $B$ as a {\em virtually central extension of $C$ by $A$}. 

\begin{proposition}\label{P1}
Let $B$ be a virtually central extension of $C$ by $A$ and  $f: G\to B$ be a quasihomomorphism whose projection to $C$ has bounded image. Then there exists a finite index subgroup $G_o<G$ such that $f|G_o$ is within finite distance from a quasihomomorphism $f_o: G\to A$ ($f_o\sim f$). 
Furthermore, if $A$ is contained in the center of $B$, then one can take $G_o=G$. 
\end{proposition}
\proof  Let $\rho: C\to Aut(A)$ denote the action of $C$ on $A$, let $Q$ be the image of $\rho$; by our assumption, the group $Q$ is finite. Without loss of generality, we may assume that the subset $f(G)$ generates  $B$ (otherwise, we 
replace $B$ with $\<f(G)\>$. By Theorem \ref{main}, we can assume (after passing to a finite index subgroup $G^o<G$ and replacing $f|G^o$ with a nearby quasihomomorphism) 
that $\Delta_f$ is contained in the center of $B$.  In particular, $\rho p(\Delta_f)=\{1\}$ and, hence,   
the composition 
$$
G \stackrel{f}{\to} B \stackrel{p}{\to} C \stackrel{\rho}{\to} Q$$ 
is a homomorphism.  Let $G_o$ denote the kernel of this homomorphism; it is a finite index subgroup of $G$. 
Then, by the construction,  $A$ is contained in the center of $B_o=f(G_o)=Ker(\rho \circ p)$.  In what follows we use the restriction of the metric from $B$ to $B_o$.

We let  $r_o: B_o\to A$ denote a nearest-point projection.  We claim that the restriction of $r_o$ 
to each $n$-neighborhood ${\mathcal N}_n(A)$ of $A$ in $B_o$ is a quasihomomorphism:
$$
r_o(xy)\sim_{S_n} r_o(x) r_o(y)
$$
for all $x, y, xy\in {\mathcal N}_n(A)$. The finite subsets $S_n$, in general, will depend on $n$. 

The proof of the claim is similar to the one in the proof of Theorem 
\ref{main}.  Let $h_i=a_i b_i\in B_o$, $a_i=r_o(h_i)$, $b_i\sim 1, b_i\in B_o, i=1,2$. Then, since $A$ is central in $B_o$, 
$$
h_1h_2= a_1 a_2 b_1 b_2,
$$
$$
r_o(h_1h_2)\sim a_1 a_2= r_o(h_1) r_o(h_2), 
$$
cf. the proof of Lemma \ref{L0}. Thus, the restriction of $r_o$ to ${\mathcal N}_n(A)$ is indeed a quasihomomorphism. Consequently, the composition $f_o=r_o\circ f: G_o\to A$ is also a quasihomomorphism. 
By the construction, the maps $f_o|G_o$ and $f|G_o$ are within finite distance from each other.  

Lastly, we note that if $A$ is central in $B$, then $Q=1$ and, thus, $B_o=B, G_o=G$. \qed 


\begin{corollary}
Suppose that $B$ is a finitely generated virtually  abelian group, $B=A\rtimes C$, where $A$ is free abelian of finite rank and $C$ is finite. Then for each quasihomomorphism $f: G\to B$, there exists a finite index subgroup $G_o<G$ such that $f|G_o$ is within finite distance from a quasihomomorphism $f_o: G\to A$. Furthermore, if $A$ is contained in the center of $B$, then one can take $G_o=G$. 
\end{corollary}

\subsection{Quasihomomorphisms to finite extensions}

Suppose that we have an extension of a group $H$, i.e., a short exact sequence 
$$
1\to K\to H \stackrel{p}{\lra} Q \to 1,
$$
and a quasihomomorphism $f: G\to H$ such that $D(f)$ is contained in the center of $H$ and $p\circ f(G)$ is finite, e.g., $Q$ is a finite group. Assume, furthermore, that the subgroup $Q_o:=p(\Delta_f)$ has finite index in 
$Q$. 

\begin{proposition}\label{P2}
Under the above assumptions, there exists a finite index subgroup $G^o< G$ and a quasihomomorphism 
$f^o: G^o\to K$, $f^o\sim f$, $D(f^o)\subset \Delta_f$.  
\end{proposition}
\proof Since the subgroup $\Delta_f$ is central in $H$, its image $Q_o=p(\Delta_f)$ is central in $Q$. The composition 
$$
G \stackrel{f}{\lra} H  \stackrel{p}{\lra} Q \to Q/Q_o 
$$
is then a homomorphism to a finite group; let $G^o$ denote its kernel. Since $p(\Delta_f)=Q_o$ and $p\circ f(G)$ is finite, 
there exists a finite subset 
$$
D_1=\{h_1, \ldots, h_n\}\subset \Delta_f,
$$
such that 
$$
f(G^o)\subset \bigcup_{i=1}^n K h_i. 
$$
Similarly to the proof of Proposition \ref{P1}, we define the retraction
$$
r: \bigcup_{i=1}^n K h_i \to K, \quad r(kh_i)=k. 
$$
Centrality of $\Delta_f$ in $H$ implies that
$$
k_1 h_{i_1} k_2 h_{i_2}= k_1 k_2 h_{i_1} h_{i_2} = k_1 k_2 h_{i_3}, 
$$
with $h_{i_1}, h_{i_2}, h_{i_3}\in D_1$. It follows that $f^o:=r\circ f|G^o$ is a quasihomomorphism and 
$$
D(f^o)\subset D_1^2 D_1^{-1}\subset \Delta_f. 
$$
Clearly, $\dist(f^o, f|G^o)<\infty$. \qed

\section{Quasihomomorphisms to hyperbolic groups}

\begin{theorem}\label{thm:hyp}
1. Suppose that $H$ is a torsion-free hyperbolic group.  Then (for an arbitrary group $G$) every unbounded quasihomomorphism $f: G\to H$  is either a homomorphism or a quasimorphism to a cyclic subgroup of $H$. 

2. Suppose that $H$ is a general hyperbolic group. Then for every unbounded quasihomomorphism $f: G\to H$ 
either the image of $f$ is contained in an elementary subgroup of $H$ or $f$ is an almost homomorphism. 
\end{theorem}
\proof In view of Corollary \ref{cor:cent}, $f(G)$ is contained in a $C$-neighborhood of the centralizer of 
$\Delta_f$ in $H$. Since $f(G)$ is infinite, it follows that the defect subgroup $\Delta=\Delta_f$ has infinite centralizer in $H$, and, hence, is elementary. Let $N=N_H(\Delta)$ be the normalizer of $\Delta$ in $H$. If 
$\Delta$ is finite then composition of $f$ with the projection to $Q=N/\Delta$ is a homomorphism and, hence, $f$ is an almost homomorphism. If $\Delta$ is infinite, then $N$ is elementary. By Lemma \ref{lem:norm}, $f(G) < N$, which concludes the proof of Part 2. 

Furthermore, suppose $H$ is torsion free. If $\Delta$ is finite, 
then it is trivial and $f$ is a homomorphism. If $\Delta$ is infinite, then $N$ is cyclic. Thus, $f: G\to N$ is a quasimorphism from $G$ to an infinite cyclic subgroup of $H$.  \qed

\medskip
The following lemma is a sharpening of the statement about quasihomomorphisms to elementary groups:

\begin{proposition}\label{prop:elementary}
If $f: G\to H$ is an unbounded quasihomomorphism to an elementary hyperbolic group $H$, then, the reduction 
$\widehat f$ of $f$ modulo the maximal finite normal subgroup $F\normal H$ either  is a quasimorphism (to $\Z$) 
or this statement holds after restricting $\widehat f$ to an index 2 subgroup $G_o<G$.  
\end{proposition}
\proof The projection of $f$, $\widehat f: G\to H/F$, is again a quasihomomorphism. Therefore,  it suffices to consider the case when $F=1$ and $H$ is either $\Z$ or $\Z_2\star \Z_2$; moreover, it suffices to consider 
the case where $H$ is generated by $f(G)$. If $H\cong \Z$, then   
$f$ is a quasimorphism. If $H\cong \Z_2\star \Z_2$, the group $\Delta_f$ has to fix the ideal boundary of $H$ pointwise (since it acts on $H$ with bounded displacement). Therefore, the composition of $f$ with 
the projection to $\Z_2$ is a homomorphism. Restricting $f$ to the kernel $G_o$ 
of this homomorphism results in a quasimorphism $\widehat f: G_o\to \Z$. \qed

\begin{corollary}
Suppose that $\Ga$ is an irreducible lattice in a semisimple Lie group of real rank $\ge 2$. 
Then every quasihomomorphism $f: \Ga\to H$, with hyperbolic target group $H$, is bounded. 
\end{corollary}
\proof First of all, it is proven in \cite{Burger-Monod} (Corollary 1.3) that $\Ga$ has only bounded quasimorphisms. Suppose, therefore, that $f: \Ga\to H$ is an unbounded  quasihomomorphism. If the image of $f$ is contained in an elementary subgroup of $H$ then, after passing to an index 2 subgroup  $\Ga_o<\Ga$, we obtain an unbounded quasimorphism $\Ga_o\to \Z$ (see Proposition \ref{prop:elementary}), which is a contradiction. 
Otherwise, by Theorem \ref{thm:hyp}, there exists a (nonelementary) subgroup 
$J<H$ such that $f(\Ga)<J$ and a finite normal subgroup $K\normal J$ such that the projection $\bar f$ of $f$ to $\bar{J}= J/K$ is a homomorphism.  The construction of quasimorphisms applied to the subgroup $J<\Ga$ (see \cite{Fujiwara}, \cite{Epstein}) yields unbounded quasimorphisms $h: J\to \Z$. Since $K$ is a normal finite subgroup in $J$, the sequence
$$
1\to K\to J\to \bar J \to 1
$$
is quasi-split (see Example \ref{ex:qsplit}) and, hence, $h$ projects to 
an unbounded  quasimorphism $\bar{h}: \bar{J}\to \Z$ (see the {\em projection} construction in \S \ref{sec:split}). Composing the  quasimorphism $\bar h$ with the homomorphism $\bar{f}: \Ga\to \bar J$, we obtain an unbounded quasimorphism $\Ga\to \Z$, which again contradicts \cite{Burger-Monod}. \qed 

\medskip
As another application of Theorem \ref{thm:hyp}, we will prove {\em deformation rigidity} of quasihomomorphisms to torsion-free hyperbolic groups. It shows that a bounded perturbation such a quasihomomorphism is seldom a quasihomomorphism. 

\begin{theorem}\label{thm:def-rig}
Suppose that $H$ is a torsion-free hyperbolic group and $f_1, f_2: G\to H$ are quasihomomorphisms with $\dist(f_1, f_2)<\infty$. Then either both $f_1, f_2$ are bounded, or both take values in the same cyclic subgroup of $H$, or $f_1=f_2$. 
\end{theorem}
\proof According to Theorem \ref{thm:hyp}, each $f_1, f_2$ is either bounded, or is a quasimorphism  to a cyclic subgroup or is a homomorphism. Recall that if $C_1, C_2$ are cyclic subgroups of a hyperbolic group $H$ then either their ideal boundaries in the Gromov boundary of $H$  are disjoint, or $C_1, C_2$ generate an elementary subgroup of $H$. In the former case, for each $R< \infty$, the intersection 
$$
{\mathcal N}_R(C_1)\cap {\mathcal N}_R(C_2)
$$
is bounded. In the setting of our theorem, it follows that if the image of $f_i$ is contained in a cyclic subgroup $C_i$ of $H$, then the image of $f_{3-i}$ is contained in a cyclic subgroup of $H$ containing $C_i$. Therefore, it remains to analyze the case when both $f_1, f_2$ are homomorphisms. 
For $x\in G$ let $C_i$ denote the cyclic subgroup of $H$ generated by $f_i(x)$. Since the homomorphisms 
$$
f_i: \< x\> \to C_i < H
$$ 
are within finite distance from each other, the subgroups $C_1, C_2$ generate a cyclic subgroup $C$ of $H$. The reader will verify that if $f_i: \< x\> \to C$ are two homomorphisms within finite distance from each other, they have to be equal. Hence, $f_1(x)=f_2(x)$ for all $x\in G$ when both $f_1, f_2$ are homomorphisms. 
\qed 

 \section{Quasihomomorphisms to $CAT(0)$ groups} 

We will need several standard facts from the theory of $CAT(0)$ groups. Recall that a group $\Ga$ is said to be a $CAT(0)$ group if there exists a $CAT(0)$ space and $\Ga\acts X$, a properly discontinuous isometric cocompact action. (This action is not required to be faithful, but the kernel of the action is necessarily finite. We are unaware, though, of any examples of $CAT(0)$ groups which do not admit faithful properly discontinuous isometric cocompact actions on $CAT(0)$ spaces.)  Recall that for an isometry $\al$ of $X$, the
{\em displacement} of $\al$ is 
$$
D_\al=\inf_{y\in X} d(y, \al  y). 
$$
Since $\Ga\acts X$ is cocompact and properly discontinuous, for every $\al\in \Ga$ this infimum is attained in $X$ and one defines the {\em minimal set} $Min_\al$ of $\al$ as 
$$
\{x\in X: d(x, \al x)= D_\al\}. 
$$
It is clear that $Min_\al$ is closed; the $CAT(0)$ property implies that $Min_\al$ is convex 
\cite[Ch II.6, Theorem 6.2]{BH} and, hence, is a $CAT(0)$ space. 


\begin{lemma}\label{lem:torsion-cat0}
Let $\Phi<\Ga$ be a finite subgroup. Then the fixed-point set $X'=X_\Phi$ of $\Phi$ in $X$ is a nonempty closed convex subspace and $X'/\Ga'$ is compact, where $\Ga'=Z_\Ga(\Phi) <\Ga$ is the centralizer of $\Phi$ in $\Ga$.  
\end{lemma}
\proof The fact that $X'$ is nonempty is a special case of the Cartan's Fixed Point Theorem (see \cite[Ch. II.2, Corollary 2.8]{BH}). The fact that $X'$ is closed is immediate; its convexity follows from uniqueness of geodesics in $CAT(0)$ spaces. Invariance of $X'$ under $\Ga'$ is again clear. 
Compactness of $X'/\Ga_{X'}$ is proven in \cite[Remark 2]{Ruane}.\qed 

\medskip 
Suppose that $\Ga\acts X$ is a cocompact properly discontinuous action of 
$\Ga$ on a $CAT(0)$ space $X$. For an abelian subgroup $A< \Ga$ consider the subset 
$Min_A\subset X$, which  is the intersection 
$$
\bigcap_{\al\in A} Min_\al.$$
Then  $Min_A$ is a nonempty closed convex subset of $X$, which splits isometrically as the direct product  $Y\times F$, where $F$ is a flat and $Y$ is a $CAT(0)$ space, see \cite[Ch II.7, Theorem 7.1]{BH}. 
Furthermore, each $\al\in A$ preserves the product decomposition of $Min_A$ and acts on $Y$ as the identity map, while $A$ acts as a cocompact group of translations on $F$. Moreover, by the same theorem, the normalizer $N_{\Ga}(A)$ of $A$ in $\Ga$ preserves $Min_A$ and contains a finite index subgroup which centralizes $A$.

\begin{lemma}\label{lem:par-cat0}
The action of $Z_{\Ga}(A)$ is cocompact on $Min_A$. In particular, the normalizer $N_{\Ga}(A)$ is finitely generated. 
\end{lemma}
\proof The subgroup $A$ is finitely generated, see e.g. \cite[Ch II.7, Corollary 7.6]{BH}. 
Let $T<A$ denote the torsion subgroup. Then $Z_T(\Ga)$ acts cocompactly on $Min_T$  by Lemma 
\ref{lem:torsion-cat0} and the group $T< A$ acts trivially on $Min_T$. Therefore, after replacing $X$ with $Min_T$, we can assume that $A$ is torsion-free. 

Now, the claim of our lemma is proven in \cite[Theorem 3.2]{Ruane} in the case when $A$ is cyclic. For a general torsion-free group $A< \Ga$, split $A$ as the product $A_1\times A_2, A_1\cong \Z$. Then the group $Z_\Ga(A_1)$ (containing $Z_\Ga(A)$) acts cocompactly on $Min_{A_1}$. The group $N_\Ga(A_1)$ preserves the product decomposition $\R\times Y_1$ of $Min_{A_1}$ and projects to a properly discontinuous cocompact group of isometries of $Y_1$. We then proceed  inductively. Since $A$ is finitely generated, lemma follows. \qed 

\begin{corollary}\label{cor:par-cat0}
The quotient $\ol{\Ga}'$ 
of $\Ga'$ by the torsion subgroup $T<A$,  contains a finite index subgroup $\Ga'_o$ 
isomorphic to $ \Pi_o \times A/T$, where $\Pi_o$ is a $CAT(0)$ group, acting properly discontinuously and cocompactly on a closed convex subset of $X$. 
\end{corollary}
\proof We have $Min_A\subset Min_T\subset X$ and the group $T$ acts trivially on $Min_T$. The centralizer of $T$ in $\Ga$ acts properly discontinuously and cocompactly on $Min_T$, hence the quotient group $\ol{\Ga}'$ also acts. The group $A/T$ is now free abelian of finite rank and   \cite[Ch II.7, Theorem 7.1]{BH} implies existence of a 
finite index subgroup $\Ga'_o$ isomorphic to $ \Pi_o \times A/T$. Lastly, the construction of the virtual splitting 
$\Ga'_o\cong \Pi_o\times A/T$ also shows that $\Pi_o$ has a properly discontinuous cocompact action on a closed convex subset of $X$: The minimal set $Min_A$ splits as a product 
$Y\times F$, where $F$ is a flat (invariant under $A$) and $Y\subset Min_A\subset X$ is a closed convex subset on which $\Pi_o$ acts properly discontinuously and cocompactly. \qed

\medskip
We can now prove our rigidity theorem for quasihomomorphisms to $CAT(0)$ groups:  

\begin{theorem}\label{thm:cat0}
Suppose that $H$ is a $CAT(0)$ group. Then for every quasihomomorphism $f: G\to H$ there exists a finite-index subgroup $G^o<G$, a $CAT(0)$ subgroup $H'<H$, a finite central subgroup $T <H'$ and a quasihomomorphism $f^o: G^o\to H'<H$ within finite distance from $f|G^o$ such that the projection 
$\bar{f}^o$ of $f^o$ to $H'/T$ splits as a product  map 
$$
{f}^o=(f_1, f_2): G^o \to  H_1 \times H_2 < H'/T,  
$$
where $f_1: G^o\to H_1$ is a homomorphism to a $CAT(0)$-group  and 
$f_2$ is a quasihomomorphism to a finitely-generated free abelian group $H_2$.  
\end{theorem}
\proof 
We continue with the notation in Theorem \ref{main1}. We obtain a finite index subgroup $G_o<G$ and a quasihomomorphism
$$
f_o: G_o\to H_o:= Z_H(\Delta_{f_o})< H
$$
within finite distance from $f|G_o$. We let $A$ denote the (finitely generated) abelian group 
$\Delta_{f_o}$ and $T< A$ the torsion subgroup. We have quotient homomorphisms
$$
H_o \stackrel{p}{\lra} H_o/T \stackrel{q}{\lra} H_o/A. 
$$
By Corollary \ref{cor:par-cat0}, $H_o/T$ contains a finite index subgroup $H^o$ 
which splits as the product $H_1\times H_2= \Pi_o\times A/T$, where $H_1=\Pi_o$ is a $CAT(0)$ group. Since $A$ contains the defect set of $f_o$, the composition $h:=q\circ p\circ f_o$  
is a homomorphism.  

Setting $H':= p^{-1}(H^o)< H_o$, we conclude that  $G^o:=h^{-1}(q(H_o)) < G_o$ 
is a finite index subgroup of $G$. Then we obtain a quasihomomorphism 
$$
f^o:= p\circ f_o =(f_1,f_2): G^o\to H_1\times H_2,
$$
where $f_1$ is a homomorphism and $f_2: G^o\to H_2$ is a quasihomomorphism to a free abelian group.  
\qed 

\begin{corollary}\label{cor:maps of lattices}
Suppose that $H$ is a uniform lattice in a connected reductive algebraic Lie group  
and $G$ is an irreducible lattice in a semisimple algebraic Lie group of real rank $\ge 2$. Then for every quasihomomorphism $f: G\to H$  there exists a finite index subgroup $G^o<G$ and a quasihomomorphism $\tilde f: G^o\to H$ within finite distance from $f|G^o$ such that $\tilde f$ is an almost homomorphism. 
\end{corollary}
\proof The group $H$ is a $CAT(0)$ group, acting (with finite kernel) on a certain nonpositively curved symmetric space. We thus can apply Theorem \ref{thm:cat0} (whose notation we will be now using). The subgroup $G^o< G$  is still an irreducible higher rank lattice; therefore, it has only bounded quasihomomorphisms to free abelian groups (see \cite{Burger-Monod}). Hence, the map $f_2$ in Theorem \ref{thm:cat0} is bounded and, therefore,
$$
\dist(f^o, f_1)< \infty, 
$$ 
$$
f_1: G^o\to H_1< H'/T
$$
is a homomorphism. Since $T$ is a finite group, the map $f_1$ lifts to an almost homomorphism 
$\tilde f: G^o\to H'< H$. By the construction, the maps $f|G^o, \tilde f$ are bounded distance from each other. \qed  

\begin{example}
There are higher rank (non-residually finite) uniform lattices $H$ as in Corollary \ref{cor:maps of lattices} with finite nontrivial center $Z_H< H$, such that $Z_H$ is contained in every finite index subgroup of $H$, see \cite{Raghunathan}. (The group $H$ is a lattice in a nonlinear connected algebraic Lie group, a $\Z_2$-central extension of the group $SO(n,2)$.)  Therefore, setting $G=H/Z_H$ and letting $f: G\to H$ be a (quasihomomorphic) lift of the identity homomorphism $G\to H/Z_H$, we obtain examples of quasihomomorphisms whose restrictions to 
any finite index subgroup $G_o<G$ are not close to homomorphisms $G_o\to H$. 
\end{example}

\begin{theorem}\label{thm:lat-rig}
Suppose that $G$ is a connected semisimple algebraic Lie group of rank $\ge 2$ without nontrivial compact normal subgroups and $\Ga< G$ is an irreducible lattice. Then  
each quasihomomorphism $f: \Ga\to \Ga$ has bounded image or is an automorphism of $\Ga$. 
\end{theorem}
\proof In view of Theorem \ref{main}, after replacing $\Ga$ with a finite index subgroup and $f$ with a nearby quasihomomorphism, we can assume that 
$$
f: \Ga\to \La< \Ga, \quad 1\to A\to \La \stackrel{p}{\lra} Q\to 1,  
$$
where $A$ is a central subgroup of $\La$, containing $\Delta_f$ and, thus, $f':=p\circ f$ is a homomorphism. 
We let $H=\bar{\La}$ denote the Zariski closure of $\La$ in $G$; we will use the notation $\bar{A}$ for the Zariski closure of $A$. 

By the Margulis Superrigidity Theorem, one of the following holds:

1. Either $f'$ has finite image, or

2. The restriction of $f'$ to a  finite index subgroup of $\Ga$ is induced by an injective homomorphism 
$G\to G'=H/\bar{A}$.  

\medskip
In the first case, the restriction of $f$ to the kernel $\Ga^o$ of $f'$ is a quasihomomorphism $\Ga^o\to A$. 
According to \cite{Burger-Monod}, $f|\Ga^o$ is bounded; hence, $f$ is bounded as well. 

In the second case, $A$ has to be finite (since $\dim(G')\ge \dim(G)$). If $A$ is nontrivial, then the dimension of $H$ is strictly smaller than the one of $G$ (since we assume that $G$ has no nontrivial normal compact subgroups). It follows that $f$ is a homomorphism in the second case.  By the Mostow Rigidity Theorem, $f$ 
is an automorphism of $\Ga$. \qed 

\section{Mapping class groups} \label{sec:mcg}

In this section we collect some definitions and facts about  mapping class groups of surfaces of finite type that will be used in the following section in order to prove a rigidity theorem for quasihomomorphisms to mapping class groups. Most of this material is quite standard, we refer the reader to \cite{FM, I} for the details. 

\subsection{Basic definitions.}

A {\em finite type} surface $\Si$ is an oriented (possibly disconnected) surface (without boundary), admitting a 
hyperbolic surface of finite area. A{\em peripheral loop} in $\Si$ is a simple loop $\al\subset \Si$ such that one of the components of $\Si \setminus \al$ is an annulus. An {\em essential multiloop} on $\Si$ is a (possibly empty) 1-dimensional compact submanifold (without boundary) $c\subset \Si$, such that no two components of $c$ are isotopic and each component of $c$ is {\em essential} in $\Si$, i.e., does not bound a disk or an annulus. If $c$ is connected and nonempty, it is called an {\em essential loop} in $\Si$. Thus, an essential loop cannot be peripheral. A subsurface $\Si'\subset \Si$ is called {\em essential} if each essential loop in $\Si'$ is still essential in $\Si$. 

\medskip 
We let $Map(\Si)$ denote the {\em mapping class group} of $\Si$, 
$$
Map(\Si)=Homeo(\Si)/Homeo_o(\Si),
$$
where $Homeo_o(\Si)$ is the connected component of the identity map $\Si\to \Si$ in the full group of homeomorphisms $Homeo(\Si)$. For $a\in Map(\Si)$ we let $h_a\in Homeo(\Si)$ denote an (unspecified) homeomorphism representing $a$.

We let $PMap(\Si)< Map(\Si)$ denote a finite index normal subgroup equal to the kernel of the homomorphism
$$
Map(\Si)\lra Aut( H_1(\Si, \Z/3)),  
$$
defined via the action of homeomorphisms of $\Si$ on its 1st homology group. We will refer to  
$PMap(\Si)$ as the {\em pure subgroup} of $Map(\Si)$, it entirely consists of {\em pure mapping classes}; we will discuss pure mapping classes in more detail in \S \ref{sec:pure}, for now we only note that each $a\in PMap(\Si)$ obviously acts trivially on $H_0(\Si)$ and preserves isotopy classes of all peripheral loops; the subgroup $PMap(\Si)$ is torsion-free. 

\medskip
Given an essential multiloop $c\subset \Si$, define the {\em twist subgroup} $T_c< PMap(\Si)$ associated to $c$, to be the group generated by Dehn twists along the components of $c$. Then $T_c$ is a free abelian group of rank $r$, where $r$ is the number of components of $c$. 

For an essential multiloop $c\subset \Si$ we let $Map_c(\Si)< Map(\Si)$ denote the subgroup consisting of mapping classes which preserves $c$ (but is allowed to permute its components and change orientation of some of the components).  The twist subgroup $T_c$ is a normal subgroup in $Map_c(\Si)$. 

If
$$
\Si=\Si_1\sqcup ... \sqcup \Si_m
$$
is a decomposition of $\Si$ into its connected components, then the group 
$Map(\Si)$ contains the product 
$$
\prod_{i=1}^m Map(\Si_i)
$$
as a finite index normal subgroup with the quotient group $Q<S_n$ (the group $Q$ acts on $\Si$ by permuting homeomorphic components of $\Si$). In the context of pure subgroups, we have
$$
PMap(\Si)\cong \prod_{i=1}^m PMap(\Si_i). 
$$

\subsection{Reduction systems and pure elements of $Map(\Si)$}\label{sec:pure} 

According to the Niel\-sen--Thurston classification, for a connected surface $\Si$ all elements of $Map(\Si)$ are classified as: 

1. Finite order. 

2. Reducible.

3. Pseudo-Anosov. 

\medskip
Each torsion subgroup of $Map(\Si)$ is finite, since the pure subgroup $PMap(\Si)$ is torsion-free.

\begin{lem}\label{lem:pA-cen}
Suppose that $\Si$ is connected. Then the normalizer $N_{Map(\Si)}(a)$ of each pseudo-Anosov element $a\in Map(\Si)$ is virtually infinite cyclic, it contains a cyclic pseudo-Anosov subgroup of finite index. The centralizer $Z_{PMap(\Si)}(a)$ of $a$ in the pure mapping class group is infinite cyclic, consisting only of pseudo-Anosov elements (and the identity). 
\end{lem}
\proof A proof can be found for instance in \cite{M}. Alternatively, the statement about centralizers in $PMap(\Si)$ is the content of  \cite[Lemma 8.13]{I}; the statement about the normalizer follows by taking the intersection
$$
Z_{PMap(\Si)}(a)= N_{Map(\Si)}(a) \cap PMap(\Si), 
$$
which has finite index in $N_{Map(\Si)}(a)$.  \qed

\begin{rem}
One also has $N_{PMap(\Si)}(a)\cong \Z$, but we will not need this property. 
\end{rem}

\begin{corollary}\label{cor:pA-cen}
Suppose that $\Si$ has the connected components $\Si_1,\ldots, \Si_m$, $a_i\in Map(\Si_i)$ are pseudo-Anosov, $i=1,\ldots,m$; define the the free abelian subgroup $A< Map(\Si)$ generated by $a_1,\ldots, a_m$. Then 
$$
Z_{PMap(\Si)}(A)\cong \Z^m. 
$$ 
\end{corollary}

Each reducible element $a\in Map(\Si)$ admits a {\em canonical reduction system} (see e.g. \cite[\S 7.4]{I}), 
which is a certain essential multiloop $c_a\subset \Si$ invariant under $h_a$ (the orientation on some of the loops can be reversed); due to the canonical nature of $c_a$,  this multiloop is invariant (up to isotopy) under the normalizer $N_{Map(\Si)}(a)$ of $a$ in $Map(\Si)$. The multiloop $c_a$ has the property that (up to isotopy) it is contained in each $h_a$-invariant multiloop in $\Si$. 

\medskip 
An element $a\in Map(\Si)$ is {\em pure} if it is orientation-preserving and either it is pseudo-Anosov or it is reducible, so that $h_a$ preserves (up to isotopy) each component of $c_a$ (together with its orientation), preserves all complementary components $\Si_i\subset \Si\setminus c_a$, and the restriction of $h_a$ to each $\Si_i$ defines either the trivial or the pseudo-Anosov element of $Map(\Si_i)$. A pure reducible element of $Map(\Si)$ is trivial iff $c_a$ is empty. Minimality of $c_a$ implies that if $a\in Map(\Si)$ is pure and preserves (up to isotopy) an essential subsurface $\Si'\subset \Si$, then $a$ preserves each component and each boundary loop of  $\Si$. The subgroup $PMap(\Si)$ consists only of pure elements, see \cite[Corollary 1.8]{I}.

\subsection{Mapping class groups of surfaces with boundary}

Suppose that $\widehat\Si$ is a surface with boundary $C$, which is a partial compactification of a finite type 
surface $\Si$, $\Si =\widehat\Si \setminus C$. In this setting one defines the {\em relative mapping class group} $Map(\widehat\Si, C)$ as the quotient, 
$$
Homeo(\widehat\Si, C)/Homeo_o(\widehat\Si, C)$$
where $Homeo(\widehat\Si, C)$  is the group of homeomorphisms of $\Si$ fixing the boundary $C$ pointwise, 
and $Homeo_o(\widehat\Si, C)< Homeo(\widehat\Si, C)$ is the identity component. We define the pure mapping class group $PMap(\widehat\Si, C)$ analogously to the case of mapping class groups for surfaces without boundary, as the kernel of the homomorphism
$$
Map(\widehat\Si, C) \to Aut(H_1(\Si, \Z/3)). 
$$

The inclusion $\Si\embed \widehat\Si$ of the interior $\Si$ of $\widehat\Si$ into $\widehat\Si$ defines the restriction homomorphism
$$
Homeo(\widehat\Si) \to Homeo(\Si)
$$
and the associated homomorphism of mapping class groups
$$
\rho: Map(\widehat\Si, C) \to Map(\Si). 
$$
The homomorphism $\rho$ is neither surjective nor injective: It's image is a finite index normal subgroup of $Map(\Si)$; the quotient $Map(\Si)/\rho(Map(\widehat\Si, C))$ is isomorphic to the semidirect product $(\Z_2)^n 
\rtimes S_n$, where $n$ is the number of the components of $C$. The kernel of $\rho$  is a free abelian 
subgroup $T_C$ of rank $n$, its free basis consists of Dehn twists $D_{\al_i}$ along loops $\al_i\subset \Si$, parallel to the components  of $C$, $i=1,\ldots, n$. However, by restricting to the pure mapping class groups we obtain a short exact sequence

\begin{equation}\label{eq:seq}
1\to T_C \to PMap(\widehat\Si, C) \to PMap(\Si)\to 1. 
\end{equation}

\begin{proposition}\label{prop:map-split}
The sequence \eqref{eq:seq} quasi-splits.
\end{proposition}
\proof The proof is by induction on the number $n$ of components of $C$.

1. Suppose that $n=1$, i.e., $C$ is connected. Let $S$ denote the surface closed surface obtained from 
$\widehat\Si$ by attaching the 2-disk along $C$. In this case, the obstruction to splitting the sequence \eqref{eq:seq} is the Euler class $e\in H^2(PMap({\Si}); \Z)$, which can be defined as the pull-back of the Euler class 
$$
\tilde{e}\in H^2(Homeo(S^1); \Z)$$ 
under the embedding 
$$
PMap(\Si) \to Aut(\pi_1(S))\to Homeo(S^1),$$
see \cite[Section 5.5.4]{FM}. The class $\tilde{e}$ is bounded, see e.g. \cite{Ghys}. Therefore, the class $e$ is bounded as well. Hence, the sequence \eqref{eq:seq} quasi-splits.

2. Suppose that the claim holds for all surfaces with $n-1$ boundary components. Let $\widehat\Si$ be a surface with 
$$
\partial \widehat\Si= C= C_1\sqcup \ldots \sqcup C_n. 
$$
Define the surface $\widehat\Si'$ by removing the circle $C_n$ from $ \widehat\Si$ and set $C':= 
C\setminus C_n= \partial \widehat\Si'$.  The surface $\widehat\Si'$ has $n-1$ boundary components, hence, 
by the induction hypothesis, there exists a quasi-splitting
$$
s':  PMap(\Si') \to PMap(\widehat\Si', C'), 
$$
of the central extension
$$
1\to T_{C'} \to PMap(\widehat\Si', C') \to PMap(\Si)\to 1. 
$$
We claim that the central extension
\begin{equation}\label{eq:1-c}
1\to T_{C_n} \to PMap(\widehat\Si, C) \to PMap(\widehat\Si', C')\to 1
\end{equation}
quasi-splits, equivalently, has bounded extension class. Given a quasi-splitting 
$$
s'':  PMap(\widehat\Si', C') \to PMap(\widehat\Si, C),
$$
we then compose it with a quasi-splitting $s'$ as above 
and obtain a quasi-splitting
$$
s= s''\circ s': PMap(\Si') \to PMap(\widehat\Si, C)
$$
of \eqref{eq:seq}. 

To prove existence of $s''$ we use the following trick. Define a new surface $S$ by attaching one-holed tori $R_1,\ldots, R_{n-1}$ to $\widehat\Si$ along each circle $C_1,\ldots, C_{n-1}$ (leaving the last circle $C_n$ untouched). The surface $S$ now has only one boundary circle. Each homeomorphism 
$$
h\in Homeo(\widehat\Si, C)
$$
extends to a homeomorphism $\tilde h$ of $S$ by the identity on each $R_i$. Projecting $\tilde{h}$ to the mapping class group $Map(S, \partial S)$,  yields embeddings 
$$
j: Map(\widehat\Si, C)\embed Map(S, C_n)  
$$
and a the analogous embedding 
$$
j: Map(\widehat\Si', C')\embed Map(S')  
$$
for the surface $S':= S\setminus C_{n}$ (which has empty boundary). 
We obtain a commutative diagram:
\begin{diagram}
1 & \rTo & T_{C_n}   &  \rTo & PMap(\widehat\Si, C)  & \rTo &   PMap(\widehat\Si', C')  & \rTo &  1\\
    &       & \dTo^{id}   &         &   \dTo^j                        &        &     \dTo^{j'}                              &        &    \\
1 & \rTo & T_{C_n}   &  \rTo & PMap(S, C_n)             & \rTo &   PMap(S')  & \rTo &  1\\   
\end{diagram}
We now apply the 1st step of induction to the bottom row of this diagram to obtain a quasi-splitting $\sigma$ of that central extension. Restricting $\sigma$ to $PMap(\widehat\Si', C') $ we obtain the desired quasi-splitting of the top row of the diagram, i.e., of the central extension \eqref{eq:1-c}. 
\qed

\subsection{Reducible subgroups}

Recall for each essential multiloop $c\subset \Si$, we have two subgroups of 
$Map(\Si)$, the subgroup $Map_c(\Si)$ and its normal subgroup $T_c$ (the twist subgroup). 
The subgroup $PMap_c(\Si):= Map_c(\Si)\cap PMap(\Si)$ still contains $T_c$.  Define the essential 
subsurface $\Si_c:= \Si \setminus c$.

\begin{lem}
The inclusion $T_c\embed PMap_c(\Si)$ defines a short exact sequence
$$
1\to T_c\to PMap_c(\Si) \stackrel{\pi}{\lra} PMap(\Si_c) \to 1. 
$$
\end{lem}
\proof The homomorphism $\pi: PMap_c(\Si) \to PMap(\Si_c)$ is induced by restricting  
homeomorphisms of $\Si$ preserving $c$ to the subsurface $\Si_c$.  The fact that its kernel contains $T_c$ 
is immediate. We next prove the equality. Let ${\mathcal N}(c)\subset \Si$ denote an open regular 
neighborhood of $c$ in $\Si$; the inclusion
$$
\Si\setminus {\mathcal N}(c) \embed \Si_c$$ 
is a homotopy-equivalence. If $f\in Homeo(\Si)$ fixes $\Si\setminus {\mathcal N}(c)$ pointwise, then $f$ projects to an element of the twist subgroup $T_c$. It follows that $\ker(\pi)= T_c$. 

To prove surjectivity of $\pi$, we note that each element of 
$$a\in PMap(\Si_c)\cong PMap(\Si\setminus {\mathcal N}(c))$$  
can be represented by a homeomorphism $h_a$ of $\Si\setminus {\mathcal N}(c)$ fixing the boundary of this subsurface pointwise. We then extend $h_a$ to each annular component of ${\mathcal N}(c)$ by an iterated Dehn twist.   The result is a homeomorphism $\tilde{h}_a$ of $\Si$ preserving $c$ are projecting to an element $\tilde{a}\in PMap_c(\Si)$ such that $\pi(\tilde{a})=a$. \qed

\subsection{Structure of infinite abelian subgroups and their normalizers} \label{sec:centralizers}

The structure of infinite abelian subgroups $A< Map(\Si)$ is described in \cite{BLM} and in \cite[chapter 8]{I}. Below is a brief review of this description. The intersection $A_P:=A\cap PMap(\Si)$ is a finite index subgroup of $A$; this subgroup is either cyclic pseudo-Anosov, or $A_P$ contains nontrivial reducible elements. We consider the latter case. For any $a_1, a_2\in A_P$, the multiloops $c_{a_1}, c_{a_2}$ are disjoint up to an isotopy, but some of the components of these multiloops could be isotopic to each other. We pick an auxiliary complete hyperbolic metric on $\Si$ and let $c_A$ denote the union of closed geodesics in $\Si$ representing 
all the loops in $c_a, a\in A_P$; then $c_A$ is an essential multiloop in $\Si$ invariant under $A_P$. Due to the canonical nature of $c_A$, this multiloop is invariant (up to isotopy) under all elements of the normalizer $N_{Map(\Si)}(A)$ of $A$ in $Map(\Si)$. In order to simplify the notation, we will denote $c_A$ by $c$. 

It follows that 
$$
T_{c}< Z_{PMap(\Si)}(A_P) < N_{Map(\Si)}(A) < Map_{c}(\Si)  
$$
and, by restricting the homomorphism $\pi$ defined in the previous section to the subgroup $N_{Map(\Si)}(A)$, 
we obtain the homomorphism 
$$
N_{Map(\Si)}(A) \stackrel{\pi}{\longrightarrow} Map(\Si_{c}) 
$$
and the exact sequence 
$$
1\to T_{c} \to N_{Map(\Si)}(A) \stackrel{\pi}{\longrightarrow} Map(\Si_{c}).  
$$

We next partition the surface $\Si \setminus c= \Si_{c}$ as 
$$
\Si_c= \Si_{c}^+\sqcup \Si_{c}^-$$ 
where $\Si_c^-$ is the union of components $\Si_i$ of $\Si_c$ such that the restriction of each $h_a, a\in A$,  to 
$\Si_i$ is isotopic to a periodic homeomorphism. In other words, a component $\Si_j$ of $\Si_c$ is contained in $\Si_c^+$ iff there exists $a\in A_P$ such that $h_a: \Si_j\to \Si_j$ is pseudo-Anosov; a component $\Si_j$ belongs 
to $\Si_c^-$ iff its stabilizer in $A$ restricts to a torsion (and, hence, finite) subgroup of $Map(\Si_j)$. 

This partition of  $\Si_c$ is preserved by $N_{Map(\Si)}(A)$ and we obtain 
 $$
 \pi=(\pi^+, \pi^-): N_{Map(\Si)}(A) {\longrightarrow} Map(\Si_c^+)\times Map(\Si_c^-) < 
 Map(\Si_c). 
 $$
Clearly, the image $\pi^\pm(N_{Map(\Si)}(A) )< Map(\Si_c^\pm)$ is contained in the normalizer
$$
N_{Map(\Si_c^\pm)}(A^\pm), \quad A^\pm=\pi^{\pm}(A). 
$$

We now restrict our discussion to pure subgroups. Set $A^\pm_P:= \pi^{\pm}(A_P)$; these are subgroups of 
$PMap(\Si_c^\pm)$. By Corollary \ref{cor:pA-cen}, the group $Z_{PMap(\Si_c^+)}(A^+_P)$ is free abelian. The fact that $A^-$ is torsion implies that $A^-P$ is trivial and, hence, $Z_{PMap(\Si_c^-)}(A^-_P)= PMap(\Si_c^-)$.  We summarize these observations as

\begin{lem}\label{lem:cen-str}
For the groups $A_P^\pm= \pi^\pm(A_P)$, we have:
$Z_{PMap(\Si_c^+)}(A^+_P)\cong \Z^r$ and $Z_{PMap(\Si_c^-)}(A^-_P)= PMap(\Si_c^-)$. Here $r=b_0(\Si_c^+)$. 
\end{lem}

\section{Quasihomomorphisms to mapping class groups}\label{sec:qh-mcg}

In this section we will extend the rigidity results from $CAT(0)$ and hyperbolic target groups to mapping class groups. The main result of this section, a rigidity theorem for quasihomomorphisms to mapping class groups is similar to Theorem \ref{thm:cat0}, except that the centralizers in mapping class groups do not (virtually) split.

\begin{theorem}\label{thm:mcg}
Suppose that $\Si$ is an oriented connected surface of finite type and $f: G\to Map(\Si)$ is a quasihomomorphism. Then there exists a finite index subgroup $G^o<G$, a quasihomomorphism 
$f^o: G^o\to Map(\Si)$, $f^o\sim f$, such that:

1. $f^o(G^o)\subset PMap_c(\Si)$ for some (possibly empty) essential multiloop $c\subset \Si$. 

2. The surface $\Si_c= \Si\setminus c$ admits a partition into subsurfaces $\Si_c= \Si_c^+\sqcup \Si_c^-$, for which we have the exact sequence
$$
1\to T_c \to PMap_c(\Si) \stackrel{(\pi^+,\pi^-)}{\lra} PMap(\Si_c^+)\times PMap(\Si_c^-)\to 1, 
$$ 
as in \S \ref{sec:centralizers}. 

3. The maps $f^\pm= \pi^\pm \circ f^o$ satisfy:

a. $f^+$ is a quasihomomorphism with free abelian target. 

b. $f^-$ is a homomorphism. 
\end{theorem}
\proof In what follows, we consider a quasihomomorphism $f: G\to Map(\Si)$ with infinite image. In view of Theorem \ref{main}, there exists a finite index subgroup $G_o<G$ and a quasihomomorphism $f_o: G_o\to Map(\Si)$, $f_o\sim f$, such that:
$$
\Delta_{f_o}< Map(\Si)
$$
is an abelian subgroup central in $\<f_o(G_o)\>$. Consider the sequence
$$
1\to PMap(\Si) \to Map(\Si) \to Aut(H_1(\Si, \Z/3))\to 1. 
$$
Applying Proposition \ref{P2} to $f_o$ and this sequence, we replace $G_o$ with its finite index subgroup $G^o$ and replace $f_o$ with a quasihomomorphism $f^o: G^o\to PMap(\Si), f^o\sim f_o$, such that 
$$
A:=\Delta_{f^o}  < \Delta_{f_o}
$$
and $f^o(G^o)$ still centralizes $A$:
$$
f^o: G^o \to Z_{PMap(\Si)}(A). 
$$
 Since the image of $f^o$ is contained in the pure mapping class group, the group $A=A_P$ is free abelian (of finite rank). If $A$ is trivial, $f^o$ is a homomorphism and we are done. 
Therefore, we will assume from now on that the group $A$ is nontrivial. 

So far, the proof is analogous to the one for CAT(0) groups. However, 
unlike in the $CAT(0)$ setting, centralizers in the mapping class group do not virtually split. 

\medskip 
There are the following possibilities for the infinite group $A$ (see \S \ref{sec:centralizers}):

{\bf 1. Pseudo-Anosov case:} There exists a pseudo-Anosov element $a\in A$. Then, 
the group $Z_{PMap(\Si)}(A)$ is infinite cyclic.  It then follows that the quasihomomorphism 
$f^o: G^o\to PMap(\Si)$ has infinite cyclic image, which concludes the proof in this case. 

\medskip 
{\bf 2. Reducible case:} $A$ contains nontrivial reducible elements. As in \S  \ref{sec:centralizers}, 
we have an $A$-invariant essential multiloop $c=c_A\subset \Si$, split the surface 
$\Si_c:=\Si\setminus c$ as $\Si_c^+\sqcup \Si_c^-$ 
and obtain  homomorphisms  
$$
Z_{PMap(\Si)}(A)< PMap_c(\Si) \stackrel{\pi}{\lra} PMap(\Si_c) = 
PMap(\Si^+_c) \times PMap(\Si^-_c),   
$$ 
$$
\pi=(\pi^+, \pi^-), \quad \pi^\pm: PMap(\Si; c)  \to PMap(\Si^\pm_c). 
$$

As we observed in Lemma \ref{lem:cen-str}, $\pi^+(Z_{PMap(\Si)}(A))\cong \Z^r$, where $r$ is the number of components of $\Si_c^+$. Therefore, for $A^+=\pi^+(A)$, we obtain the quasihomomorphism 
$$
f^+=\pi^+\circ f^o: G^o \to Z_{PMap(\Si^+_c)}(A^+) \cong \Z^r. 
$$
As for $\Si_c^-$, the projection $\pi^-(A)$ is trivial and, since $A$ contains the defect subgroup of $f^o$, the composition
$$
f^-=\pi^-\circ f^o: G^o \to Z_{PMap(\Si^-_c)}(A^-) = PMap(\Si^-_c) 
$$
is a homomorphism. \qed

\begin{corollary}\label{cor:mcg}
Suppose that $\Gamma$ is an irreducible lattice in a connected semisimple Lie group of rank $\ge 2$, without compact factors. 
Then every quasihomomorphism of $\Gamma$ to a mapping class group $Map(\Si)$ has finite image. 
\end{corollary}
\proof Suppose to the contrary that $f: \Ga\to Map(\Si)$ is an unbounded quasihomomorphism. As in 
Theorem \ref{thm:mcg}, we replace $\Ga$ with its finite index subgroup $\Ga^o$ (which is still an irreducible lattice of rank $\ge 2$) and replace $f$ with $f^o\sim f, f^o: \Ga^o\to PMap(\Si)$. The compositions 
$$
f^\pm=\pi^\pm \circ f^o: \Ga^o\to PMap(\Si_c^\pm), 
$$ 
satisfy the property that $f^+$ is a quasihomomorphism to a free abelian group $A_1$ 
and $f^-$ is a homomorphism. 
The homomorphism $f^-$ has to have finite image (see \cite{BF, FMas, KM}); actually, in our setting, the image of $f^-$ is trivial since $PMap(\Si_c^-)$ is torsion-free. Therefore, the image of the map $f^o$ is contained in the abelian subgroup $B< PMap(\Si)$, the preimage $(\pi^{+})^{-1}(A_1)$. Therefore, $f^o$ is  bounded in view of \cite{Burger-Monod}. A contradiction. \qed

\section{Quasihomomorphisms to groups acting trees}

Suppose $T$ is a simplicial tree and $H=Aut(T)$ is the group 
of automorphisms of $T$ acting on $T$ without inversions. 

\begin{definition} 
Suppose that $T'\subset T$ is a nonempty simplicial subtree and that $f: G\to Aut(T)$ is 
a quasi-homomorphism whose image preserves $T'$. Let $H'= Aut_{T'}(T)$ denote the subgroup of $Aut(T)$ preserving $T'$. We have the restriction homomorphism $r:  H'\to Aut(T')$. The composition 
$f':= r\circ f$ is a quasihomomorphism $f': G\to  Aut(T')$. In this situation we will say that the quasihomomorphism 
$f$ is a {\em lift} of the quasihomomorphism $f'$. 
\end{definition} 

We now proceed with the analysis of quasihomomorphisms $f: G\to  H=Aut(T)$. 
Using Theorem \ref{main1}, we find  $f_o: G_o \to H_o =\langle f_o(G_o) \rangle$, such that  
$\Delta=\Delta_{f_o}$ is central in $H_o$.

{\bf Case 1. Axial case:} Suppose that $\Delta$ contains an axial isometry $\delta$ of $T$, i.e., an isometry 
which preserves a complete geodesic $T'$ in $T$ and acts on $T'$ as a nontrivial translation, i.e., $T'$ 
is the {\em axis} of $\delta$. Since each axial isometry has unique axis, the axis $T'$ 
of $\delta$ is invariant under $H_o$ and $H_o$ acts on $L$ by integer translations. (Centrality of  $\Delta$ implies that every element of $H_o$ preserves the orientation on $T'$.) Let 
$$
Aut_{T'}^+(T)< Aut_{T'}(T)$$ 
denote the subgroup of $Aut(T)$ preserving $T$ and its orientation. We have a natural homomorphism 
$$
\tau: Aut^+_{T'}(T) \to \Z,
$$
sending each $h\in Aut^+_{T'}(T)$  to the translation number for its action on $T'$. 
Composing $f_o$ with $\tau$ we obtain a quasimorphism 
$$
f_o'=\tau\circ f_o: G_o \to \Z.$$ 
Thus, in this setting, $f_o$ is a lift of a quasihomomorphism to $\Z$.  

\medskip 
{\bf Case 2. Elliptic case:} Suppose that  $\Delta$ contains only elliptic isometries, i.e., each element of $\Delta$ has a fixed point in $T$. Recall that the defect group $\Delta$ is finitely generated abelian. 
 
\begin{lemma}
Let $A$ be a finitely generated abelian group acting isometrically on a tree $T$ such that every element of $A$ is elliptic. Then the fixed-point set of the action of $A$ on $T$ is nonempty. 
\end{lemma} 
 \proof We let $A_1,\ldots, A_n$ denote cyclic factors of $A$. The fixed subtree $T_i$ of each $A_i$ is nonempty. We claim that the tree 
 $$
 T'= T_1\cap \ldots \cap T_n
 $$
is nonempty. The proof is by induction on $n$. The claim is clear for $n=1$. Assume that it holds for $n-1$. 
The subgroup $A'< A_1\times \ldots \times A_{n-1}< A$ preserves the tree $T_n$ and each element  of $A'$ acts on $T_n$ as an elliptic isometry. Thus, the claim follows from the induction hypothesis. \qed 
 
 \medskip 
Applying this lemma to the group  $A=\Delta_{f_o}$, we conclude that its fixed-point set in $T$ 
is a nonempty subtree $T'\subset T$. By the normalization property, this subtree has to be invariant under 
$H_o$ and, as above, we obtain the homomorphism 
$$
 f_o' =r\circ f_o: G_o\to H'=Aut(T').
$$  
Hence, the quasihomomorphism $f_o$ is a lift of the homomorphism $f_o'$. 

This proves 

\begin{lemma}\label{lem:tree}
If $f: G\to H=Aut(T)$ is a quasi-homomorphism then, there exists $f_o: G_o\to H$, $f_o\sim f$, 
such that: 

1. Either $f_o$ is an extension of a quasimorphism  $f'_o: G_o\to \Z< H$,  or 

2. $f_o$ is an extension of a homomorphism $f'_o: G_o\to H'=Aut(T')$ where $T'\subset T$ is a nonempty subtree. 
\end{lemma}

\begin{corollary}\label{cor:fixedvertex}
Suppose that $G_o$ has no unbounded quasimorphisms and satisfies the property FA (e.g., $G$ 
is an irreducible lattice in a connected semisimple Lie group of rank $\ge 2$). Then there exists a subgroup $G^o<G_o$ of finite index and a quasihomomorphism $f^o: G^o\to Aut(T), f^o\sim f_o$, 
such that $f^o(G^o)$ fixes a vertex in $T$. 
\end{corollary}
\proof Since $G_o$ satisfies the property FA, $f_o(G)$ has a fixed vertex in $T'$ in the {\em elliptic case}. 
Hence, in this situation, we can take $G^o=G_o, f^o=f_o$. Consider now the hyperbolic case. 
By the assumptions, the quasimorphism $f_o': G_o\to \Z$ has finite image. Therefore, we apply 
Proposition \ref{P2} to the exact sequence
$$
1\to K \to Aut^+_{T'}(T) \stackrel{\tau}{\lra} \Z\to 1
$$
and conclude that there exists a finite index subgroup $G^o< G_o$ and a quasihomomorphism 
$f^o: G^o\to K$ with $f^o\sim f_o$. The image of $f^o$ fixes each vertex of $T'$. \qed 

\begin{corollary}
Suppose that $H$ is the fundamental group of a graph of groups where every vertex group is hyperbolic. 
Then for every group $G$ satisfying the hypothesis of Corollary \ref{cor:fixedvertex}, each quasihomomorphism 
$f: G\to H$ has finite image. 
\end{corollary}

\section{Other generalizations of homomorphisms}\label{sec:HS}

In this section we compare the notion of quasihomomorphisms used in this paper and going back to Ulam, with several other notions. In order to avoid the notation confusion, we will refer to quasihomomorphisms used earlier as {\em Ulam--quasi\-homomorphisms}. The  other notions discussed in this section are equivalent to the one of Ulam--quasi\-homomorphism when the target is $\Z$, but differ in general.  

\subsection{Algebraic and geometric quasihomomorphisms} 

Let $G$ and $H$ be groups and $d$ is a left-invariant metric on $H$. A map $f: G\to H$ is an {\em algebraic quasihomomorphism} if there exists  a bounded subset  $S\subset H$ such that for all $x, y\in G$ we have:
$$
f(xy)= s_1 f(x) s_2 f(y) s_3, \quad s_i\in S, i=1, 2, 3. 
$$
The true novelty in this definition (comparing to the one of Ulam--quasi\-homo\-mor\-phisms) is presence of the element $s_2$. This class of maps is preserved by the following {\em bi-bounded perturbation} procedure: Pick a bounded subset $B\subset (H,d)$ and consider a map $f': G\to H$ such that for each $x\in G$, $f(x)\in B f(x) B$. Then $f'$ is again an 
algebraic quasihomomorphism.

Alternatively, one can require the more restrictive condition
$$
f(xy)= f(x) s_2 f(y) s_3, \quad s_i\in S, i=2, 3, 
$$
where $S$ is a bounded subset of $(H,d)$. We refer to such maps as {\em geometric quasihomomorphisms}. 
Geometric and algebraic quasihomomorphisms are stable under bounded perturbations. 
This presents a sharp contrast with Ulam's quasihomomorphisms (cf. Theorem \ref{thm:def-rig}). 

We let $AQHom(G, (H,d))$ and $GQHom(G, (H,d))$ denote the sets of algebraic and geometric  quasihomomorphisms, and denote by $UQHom(G, (H,d))$ the set of Ulam-quasihomomorphisms. 

\begin{example}
1. Each map $f: H\to H$ such that $\dist(f, id)<\infty$, is a geometric quasihomomorphism. 

2. Compositions of algebraic (respectively, geometric) quasihomomorphisms are again (respectively, geometric)  quasihomomorphisms. 
\end{example}


We will give some interesting examples of geometric quasihomomorphisms in the next section. 

A  situation when geometric quasihomomorphisms appear naturally is the one of Margulis-type superrigidity: Suppose that $\Ga< G$ is a uniform lattice in a connected Lie group (equipped with a left-invariant Riemannian metric) and $\phi: \Ga\to (H,d)$ is a homomorphisms. Then for a nearest-point projection $\rho: G\to \Ga$ (which is a geometric quasihomomorphism), the composition 
$$
f= \phi\circ \rho: G\to (H,d)
$$
is again a geometric quasihomomorphism. If $G$ is a simple noncompact group of rank $\ge 2$, then  the Margulis Superrigidity Theorem implies that such geometric quasihomomorphism $f$ is within finite distance from a homomorphism $G\to H$, provided that $H$ is another connected Lie group (and $d$ is induced by a left-invariant Riemannian metric on $H$). This leads to:

\begin{question}\label{ques:a-rig}
Suppose that $G$ is a connected simple Lie group of real rank $\ge 2$ and $(H,d)$ is a connected Lie group with trivial center, equipped with a metric $d$ induced by a left-invariant Riemannian metric on $H$. 
Is it true that {\em every} geometric quasihomomorphism $f: G\to (H,d)$ is within finite distance from a homomorphism?  
\end{question}

Note that the answer is clearly negative for all rank 1 Lie groups, for instance, because these groups contain uniform lattices admitting unbounded quasimorphisms to $\Z$.  

\begin{prob}
Describe $AQHom(G, H)$ for simple connected Lie groups $G, H$ of rank $\ge 2$. Is it true that 
each $f\in AQHom(G,G)$  is a bi-bounded perturbation of a homomorphism? 
\end{prob}

\subsection{Middle--quasihomomorphisms}

The following definition is inspired by a correspondence  from Narutaka Ozawa. 

\begin{definition}
A map $f: G\to H$ of two groups is a {\em middle--quasihomomorphism} if there exists a finite subset $S\subset H$ such that for all $x, y\in G$, there is $s\in S$ satisfying
$$
f(xy)= f(x) s f(y). 
$$ 
We let $MQHom(G, H)$ denote the set of all middle--quasihomomorphisms $G\to H$. 
\end{definition}

By the definition, each middle--quasihomomorphism is geometric. As with other  quasihomomorphisms, composition preserves middle--quasiho\-momor\-phisms.

\medskip
Below is an interesting construction of {\em middle--quasihomomorphisms} $f: F_2\to F_2$ 
which is a generalization of the Brooks' construction of quasimorphisms of free groups. 
Let $a, b$ be  free generators of the free group $F_2$. We say that two nonempty reduced words $u, v$ in the 
alphabet $a^{\pm 1}, b^{\pm 1}$ are {\em non-overlapping} if for every reduced word $w$ 
in the alphabet $a^{\pm 1}, b^{\pm 1}$ any two subwords which are copies of distinct elements of 
 $$
T=\{u, u^{-1}, v, v^{-1}\}, 
$$
are disjoint. (In particular,  the words $u, v$ are cyclically reduced.) For instance, for $m\ge 2$  the words 
$$
u=a^m b a^m, \quad v= b^m a b^m 
$$
satisfy this condition. Let $L$ denote the maximum of lengths of $u$ and $v$. 

The subgroup $H$ generated by $u$ and $v$ is free of rank 2 (with the generators $u, v$), since this subgroup cannot be cyclic. 

Now, given a reduced word $w$  in the alphabet $a^{\pm 1}, b^{\pm 1}$, consider all the subwords 
$t_1,\ldots, t_n$ (listed in the order of their appearance in $w$) which belong to the set $T$. Define the map
$$
f: F_2\to H,
$$
$$
f(w)=f_{u,v}(w):= t_1 \ldots t_n\in F_2.  
$$  
If $n=0$, we set $f(w)=1$. Let $\al: H\to \Z$ denote the homomorphism 
sending $v$ to $0\in \Z$ and $u$ to $1\in \Z$. Then the composition $\be=\al\circ f$  is  the Brooks quasimorphism $F_2\to \Z$, associated with the word $u$. 

\medskip 
It is clear from the construction that $f(w)= (f(w))^{-1}$ for each $w\in F_2$.

\begin{theorem} 
1. $f$ is a middle--quasihomomorphism. 

2. The image of $f$ is infinite and is not contained in the $R$-neighborhood of an infinite cyclic subgroup of $F_2$ for any $R<\infty$.

3. The map $f$ is not within finite distance from a homomorphism. 
\end{theorem}
\proof 1. We first check that $f$ is a middle--quasihomomorphism. Consider two words $w_1, w_2$, $w_1= w_1' w_1'', w_2= w_2' w_2''$, where $w_1'', w_2'$ are maximal with the property that 
$$
w_1'' w_2'=1. 
$$ 
We let $J(w_i)$ denote the ordered set (listed in the order of their appearance in $w_i$) of
 subwords in $w_i$ which are copies of elements of $T$   overlapping both $w_i', w_i''$. 
 Then the ordered product $Y_i$ of the elements of $J(w_i)$ has length 
 $\le L^2$. Furthermore, 
$$
f(w_1)= X_1 Y_1 Z_1, \quad f(w_2)= Z_1^{-1} Y_2 Z_2, 
$$
and for the element $w_3\in F_2$ represented by $w_1 w_2$ we have
$$
f(w_3) = X_1 Y_3 Z_2,
$$
where $|Y_3|\le L^2$. Set 
$$
s_2= Y_1^{-1} Y_3 Y_2^{-1}. 
$$ 
Then
$$
f(w_3)= f(w_1) s_2 f(w_2), 
$$
where $s_2$ has length $\le 3L^2$. This proves the first claim. 

\medskip 
2. It is clear that $f(u^n)= u^n$ and $f(v^n)= v^n$ for each  $n$. Since the cyclic subgroups of $F_2$ generated by $u$ and by $v$ are not  Hausdorff-close, the second claim of the theorem follows. 

3. The map $f$ sends both cyclic subgroups $\<a\>$ and $\<b\>$ to $\{1\}$. Therefore, for each map 
$f': F_2\to F_2$ within finite distance from $f$, the images of $\<a\>$ and $\<b\>$ are bounded. Hence, $f'$ can be a homomorphism only if it is the constant map $F_2\to \{1\}$. Since $f$ is unbounded, we conclude that it cannot be within  finite distance from a homomorphism. \qed 


\subsection{Quasimorphisms of Hartnick and Schweitzer}

In their paper \cite{HS}, which appeared shortly after the initial version of our paper was posted, 
Hartnick and Schweitzer introduce the following notion, which we will refer to as an HS--quasimor\-phism:

\begin{definition}\label{def:HS} 
A map $f: G\to H$ of two groups is an HS--quasimorphism if for each quasimorphism $\varphi: H\to \R$, 
the composition $\varphi\circ f: G\to \R$ is a quasimorphism. (Note that $H$ need not be equipped with a metric.) 
We will use the notation $HSQMor(G, H)$ for the set of HS--quasimorphisms. 
\end{definition}

In other words, Hartnick and Schweitzer take the concept of quasimorphisms (quasihomomorphisms to $\R$) as central, and then define HS--quasimorphisms in a categorical fashion. It is immediate that composition preserves  HS--quasihomomorphisms.  If we equip the target group $H$ with a discrete proper left-invariant metric (whose  choice is irrelevant and will be suppressed), then, clearly,  
$$
UQHom(G, H) \subset GQQHom(G, H) \subset AQHom(G, H) \subset HSQMor(G, H), 
$$
$$
MQHom(G, H) \subset GQQHom(G, H) \subset AQHom(G, H) \subset HSQMor(G, H).  
$$
In particular, as with algebraic quasihomomorphisms,  if $f_1: G\to H$ is an HS--quasihomomorphism and 
$\dist(f_1, f_2)<\infty$, then $f_2: G\to H$ is again an HS--quasihomomorphism. Hartnick and Schweitzer prove, among other interesting results,  that free groups $F_n$ of finite rank $n\ge 2$ have abundant supply of HS--automorphisms. More precisely, let $QAut(F_n)$ denote the space of HS--quasiautomorphism $F_n\to F_n$, 
$Hom(F_n, \R)$ is the space usual homomorphisms and ${\mathcal H}(F_n)$ the space of homogeneous quasimorphisms $F_n\to \R$. Then, according to Theorem 1 of \cite{HS}, the closure of the linear span 
of the $QAut(F_n)$-orbit of  $Hom(F_n, \R)$ is the entire space ${\mathcal H}(F_n)$. 

A drawback of  Definition \ref{def:HS} is that it is only meaningful for maps to groups $H$ which admit abundant supply of quasimorphisms (e.g., hyperbolic groups). In contrast, if $H$ is an irreducible lattice of  rank $\ge 2$, then {\em every map} $G\to H$ is an HS--quasimorphism, as $H$ has only bounded quasimorphisms. 
In contrast, Theorem \ref{thm:lat-rig} shows that if $\Ga< G$ is an irreducible lattice in a  connected semisimple Lie group $G$ of rank $\ge 2$, without nontrivial compact normal subgroups, then each Ulam-quasihomomorphism $f: \Ga\to \Ga$ has finite image or is an automorphism.

Addresses:

\noindent K.F. : Department of Mathematics, \\
Kyoto University, Kyoto, 606-8502 Japan \\
email : kfujiwara@math.kyoto-u.ac.jp
\\

\noindent M.K.: Department of Mathematics, \\
University of California, Davis\\
CA 95616, USA\\
email: kapovich@math.ucdavis.edu

\end{document}